\documentclass[11pt]{tran-l}
\usepackage[ansinew]{inputenc}
\usepackage[centertags]{amsmath}
\usepackage{exscale}
\usepackage{relsize}
\usepackage{amsfonts}
\usepackage{amssymb}
\usepackage{amsthm}
\usepackage{newlfont}

\mathchardef\ordinarycolon\mathcode`\:
  \mathcode`\:=\string"8000
  \begingroup \catcode`\:=\active
    \gdef:{\mathrel{\mathop\ordinarycolon}}
  \endgroup

\def\Z{{\Bbb Z}}

\def\R{{\Bbb R}}
\newtheorem*{mthm}{Main Theorem}
\newtheorem{thm}{Theorem}
\newtheorem{lem}[thm]{Lemma}
\newtheorem{prop}[thm]{Proposition}

\newtheorem{rem}[thm]{Remark}
\newtheorem{Cor}[thm]{Corollary}

\newcommand{\D}{\mathcal{D}}
\newcommand{\B}{\mathcal{B}}
\newcommand{\g}{\mathcal{G}}
\newcommand{\h}{\mathcal{H}}
\newcommand{\aut}{\emph{Aut}}
\newcommand{\Aut}{\mbox{Aut}}
\newcommand{\Out}{\mbox{Out}}

\newcommand{\Si}{\mathit{\Sigma}}
\newcommand{\q}{\overline{q}}

\begin{document}
\title[Flag-transitive Steiner \mbox{4-designs}]
{The Classification of Flag-transitive\\ Steiner \mbox{4-Designs}}

\author{Michael Huber}

\address{Mathematisches Institut der Universit\"{a}t T\"{u}bingen, Auf der Morgenstelle~10,
D-72076~T\"{u}bingen, Germany}

\email{michael.huber@uni-tuebingen.de}

\subjclass[2000]{Primary 51E10; Secondary 05B05, 20B25}

\keywords{Steiner designs, flag-transitive group of automorphisms,
\mbox{$2$-transitive} permutation group}

\thanks{The author gratefully acknowledges partial support by the Deutsche Forschungsgemeinschaft (DFG)}

\date{March 15, 2006; and in revised form November 23, 2006}

\commby{H. Van Maldeghem}


\begin{abstract}
Among the properties of homogeneity of incidence structures
flag-transitivity obviously is a particularly important and natural
one. Consequently, in the last decades flag-transitive Steiner
\mbox{$t$-designs} (i.e. flag-transitive $t$-$(v,k,1)$ designs) have
been investigated, whereas only by the use of the classification of
the finite simple groups  has it been possible in recent years to
essentially characterize all flag-transitive Steiner
\mbox{$2$-designs}. However, despite the finite simple group
classification, for Steiner \mbox{$t$-designs} with parameters $t>2$
such characterizations have remained challenging open problems for
about 40 years (cf.~\cite[p.\,147]{Del1992}
and~\cite[p.\,273]{Del1995}, but presumably dating back to around
1965). The object of the present paper is to give a complete
classification of all flag-transitive Steiner \mbox{$4$-designs}.
Our result relies on the classification of the finite doubly
transitive permutation groups and is a continuation of the author's
work~\cite{Hu2001,Hu2005} on the classification of all
flag-transitive Steiner \mbox{$3$-designs}.
\end{abstract}

\maketitle

\section{Introduction}\label{intro}

For positive integers $t \leq k \leq v$ and $\lambda$, we define a
\mbox{\emph{$t$-$(v,k,\lambda)$ design}} to be a finite incidence
structure \mbox{$\D=(X,\B,I)$}, where $X$ denotes a set of
\emph{points}, $\left| X \right| =v$, and $\B$ a set of
\emph{blocks}, $\left| \B \right| =b$, with the properties that each
block $B \in \B$ is incident with $k$ points, and each
\mbox{$t$-subset} of $X$ is incident with $\lambda$ blocks. A
\emph{flag} of $\D$ is an incident point-block pair, that is $x \in
X$ and $B \in \B$ such that $(x,B) \in I$. We consider automorphisms
of $\D$ as pairs of permutations on $X$ and $\B$ which preserve
incidence, and call a group \mbox{$G \leq \Aut (\D)$} of
automorphisms of $\D$ \emph{flag-transitive} (respectively
\emph{block-transitive}, \emph{point $t$-transitive}) if $G$ acts
transitively on the flags (respectively transitively on the blocks,
$t$-transitively on the points) of $\D$. For short, $\D$ is said to
be, e.g., flag-transitive if $\D$ admits a flag-transitive group of
automorphisms.

For historical reasons, a \mbox{$t$-$(v,k,\lambda)$ design} with
$\lambda =1$ is called a \emph{Steiner \mbox{$t$-design}} (sometimes
this is also known as a \emph{Steiner system}). We note that in this
case each block is determined by the set of points which are
incident with it, and thus can be identified with a $k$-subset of
$X$ in a unique way. If $t<k<v$ holds, then we speak of a
\emph{non-trivial} Steiner \mbox{$t$-design}.

Among the properties of homogeneity of incidence structures
flag-\linebreak transitivity obviously is a particularly important
and natural one. Consequently, in the last decades flag-transitive
Steiner \mbox{$t$-designs} have been investigated, in particular for
the case $t=2$. The general study of flag-transitive Steiner
\mbox{$2$-designs} was introduced by D.~G.~Higman and
J.~E.~McLaughlin~\cite{HigMcL1961} proving that a flag-transitive
group $G \leq \Aut(\D)$ of automorphisms of a Steiner
\mbox{$2$-design} $\D$ is necessarily primitive on the points of
$\D$. They posed the problem of classifying all finite
flag-transitive projective planes, and showed that such planes are
desarguesian if its orders are suitably restricted. Much later
W.~M.~Kantor~\cite{Kant1987} determined all such planes apart from
the still open case when the group of automorphisms is a Frobenius
group of prime degree
(cf.~\cite[Sect.\,1,\,2]{Buek1988},~\cite{Kant1985b}, and
~\cite{Kant1993} for a detailed survey on flag-transitive planes).
In a big common effort, F.~Buekenhout, A.~Delandtsheer, J.~Doyen,
P.~B.~Kleidman, M.~W.~Liebeck, and
J.~Saxl~\cite{Buek1990,Del2001,Kleid1990,Lieb1998,Saxl2002}
essentially characterized all finite flag-transitive linear spaces,
that is flag-transitive Steiner \mbox{$2$-designs} (for the
incomplete case with a $1$-dimensional affine group of
automorphisms, see~\cite[Sect.\,4]{Buek1990}
and~\cite[Sect.\,3]{Kant1993}). Their result, which was announced in
1990, makes use of the classification of the finite simple groups.

However, despite the classification of the finite simple groups, for
Steiner \mbox{$t$-designs} with parameters $t>2$ such
characterizations have remained challenging open problems for about
40 years (see~\cite[p.\,147]{Del1992} and~\cite[p.\,273]{Del1995},
but presumably dating back to around 1965,
cf.~\cite{Luene1965b,Tits1964}). Recently, the
author~\cite{Hu2001,Hu2005} completely determined all
flag-transitive Steiner \mbox{$3$-designs} using the classification
of the finite doubly transitive permutation groups, which in turn
relies on the classification of the finite simple groups.

The object of the present paper is to give a complete classification
of all flag-transitive Steiner \mbox{$4$-designs}. Our approach uses
again the classification of the finite doubly transitive permutation
groups.

\medskip

The classification of all non-trivial Steiner \mbox{$4$-designs}
admitting a flag-transitive group of automorphisms can be stated as
follows.

\begin{mthm}
Let $\D=(X,\B,I)$ be a non-trivial Steiner \mbox{$4$-design}. Then
\mbox{$G \leq \aut(\D)$} acts flag-transitively on $\D$ if and only
if one of the following occurs:

\begin{enumerate}
\item[(1)] $\D$ is isomorphic to the Witt \mbox{$4$-$(11,5,1)$}
design, and \mbox{$G \cong M_{11}$},

\medskip

\item[(2)] $\D$ is isomorphic to the Witt \mbox{$4$-$(23,7,1)$}
design, and \mbox{$G \cong M_{23}$}.
\end{enumerate}
\end{mthm}

\medskip

For a detailed description of the \emph{Witt} \mbox{$t$-$(v,k,1)$}
\emph{designs} with their associated \emph{Mathieu groups} $M_v$ of
degree $v$, we refer, e.g., to~\cite{Witt1938}.

\bigskip


\section{Definitions and Preliminary Results}\label{Prelim}

If $\D=(X,\B,I)$ is a \mbox{$t$-$(v,k,\lambda)$} design with $t \geq
2$, and $x \in X$ arbitrary, then the \emph{derived} design with
respect to $x$ is \mbox{$\D_x=(X_x,\B_x, I_x)$}, where $X_x = X
\backslash \{x\}$, \mbox{$\B_x=\{B \in \B: (x,B)\in I\}$} and $I_x=
I \!\!\mid _{X_x \times \; \B_x}$. In this case, $\D$ is also called
an \emph{extension} of $\D_x$. Obviously, $\D_x$ is a
\mbox{$(t-1)$-$(v-1,k-1,\lambda)$} design.

Let $G$ be a permutation group on a non-empty set $X$. For $g \in
G$, let $\mbox{Fix}_X(g)$ denote the set of fixed points of $g$ in
$X$. We call $G$ \emph{semi-regular} if the identity is the only
element that fixes any point of $X$. If additionally $G$ is
transitive, then it is said to be \emph{regular}. Furthermore, for
$x \in X$, the orbit $x^G$ containing $x$ is called \emph{regular}
if it has length $\left| G \right|$. If \mbox{$\{x_1,\ldots,x_m\}
\subseteq X$}, let $G_{\{x_1,\ldots,x_m\}}$ be its setwise
stabilizer and $G_{x_1,\ldots,x_m}$ its pointwise stabilizer.

For \mbox{$\D=(X,\B,I)$} a Steiner \mbox{$t$-design} with \mbox{$G
\leq \Aut (\D)$}, let $G_B$ denote the setwise stabilizer of a block
$B \in \B$, and for $x \in X$, we define $G_{xB}= G_x \cap G_B$.

Let $\Z^+$ be the set of positive integers (without $0$). For
integers $m$ and $n$, let $(m,n)$ denote the greatest common divisor
of $m$ and $n$, and we write $m \mid n$ if $m$ divides $n$.

For any $x \in \R$, let $\lfloor x \rfloor$ denote the greatest
positive integer which is at most $x$.

All other notation is standard.

\medskip

The approach to the classification of all flag-transitive Steiner
\mbox{$4$-designs} starts with the following proposition which can
be deduced from a result of Block~\cite[Thm.\,2]{Block1965}.

\begin{prop}{\em (cf.~\cite{Buek1968,Hu2005}).}\label{flag2trs}
Let $\D=(X,\B,I)$ be a Steiner \mbox{$t$-design} with $t \geq 3$.
If \mbox{$G \leq \aut(\D)$} acts flag-transitively on $\D$, then
$G$ also acts point \mbox{$2$-transitively} on $\D$.
\end{prop}

We note that if $t=2$, then it is elementary that conversely the
point $2$-transitivity of \mbox{$G \leq \Aut(\D)$} implies its
flag-transitivity.

The above result allows us to make use of the classification of all
finite \mbox{$2$-transitive} permutation groups, which itself relies
on the classification of all finite simple groups
(cf.~\cite{CSK1976,Gor1982,Her1974,Her1985,Hup1957,Kant1985,Lieb1987,Mail1895}).

The list of groups is as follows.

Let $G$ be a finite \mbox{$2$-transitive} permutation group on a
non-empty set $X$. Then $G$ is either of

{\bf (A) Affine Type:} $G$ contains a regular normal subgroup $T$
which is elementary Abelian of order $v=p^d$, where $p$ is a
prime. \mbox{If $a$ divides $d$,} and if we identify $G$ with a
group of affine transformations
\[x \mapsto x^g+u\]
of $V=V(d,p)$, where $g \in G_0$ and $u \in V$, then particularly
one of the following occurs:

\pagebreak

\begin{enumerate}

\smallskip

\item[(1)] $G \leq A \mathit{\Gamma} L(1,p^d)$

\smallskip

\item[(2)] $G_0 \unrhd SL(\frac{d}{a},p^a)$, $d \geq 2a$

\smallskip

\item[(3)] $G_0 \unrhd Sp(\frac{2d}{a},p^a)$, $d \geq 2a$

\smallskip

\item[(4)] $G_0 \unrhd G_2(2^a)'$, $d=6a$

\smallskip

\item[(5)] $G_0 \cong A_6$ or $A_7$, $v=2^4$

\smallskip

\item[(6)] $G_0 \unrhd SL(2,3)$ or $SL(2,5)$, $v=p^2$,
$p=5,7,11,19,23,29$ or $59$, or $v=3^4$

\smallskip

\item[(7)] $G_0$ contains a normal extraspecial subgroup $E$ of
order $2^5$, and $G_0/E$ is isomorphic to a subgroup of $S_5$,
$v=3^4$

\smallskip

\item[(8)] $G_0 \cong SL(2,13)$, $v=3^6,$
\end{enumerate}

\smallskip

or

\medskip

{\bf  (B) Almost Simple Type:} $G$ contains a simple normal subgroup
$N$, and \mbox{$N \leq G \leq \Aut(N)$}. In particular, one of the
following holds, where $N$ and $v=|X|$ are given as follows:
\begin{enumerate}

\smallskip

\item[(1)] $A_v$, $v \geq 5$

\smallskip

\item[(2)] $PSL(d,q)$, $d \geq 2$, $v=\frac{q^d-1}{q-1}$, where
$(d,q) \not= (2,2),(2,3)$

\smallskip

\item[(3)] $PSU(3,q^2)$, $v=q^3+1$, $q>2$

\smallskip

\item[(4)] $Sz(q)$, $v=q^2+1$, $q=2^{2e+1}>2$ \hfill (Suzuki
groups)

\smallskip

\item[(5)] $Re(q)$, $v=q^3+1$, $q=3^{2e+1} > 3$ \hfill (Ree
groups)

\smallskip

\item[(6)] $Sp(2d,2)$, $d \geq 3$, $v = 2^{2d-1} \pm 2^{d-1}$

\smallskip

\item[(7)] $PSL(2,11)$, $v=11$

\smallskip

\item[(8)] $PSL(2,8)$, $v=28$ (N is not \mbox{$2$-transitive)}

\smallskip

\item[(9)] $M_v$, $v=11,12,22,23,24$ \hfill (Mathieu groups)

\smallskip

\item[(10)] $M_{11}$, $v=12$

\smallskip

\item[(11)] $A_7$, $v=15$

\smallskip

\item[(12)] $HS$, $v=176$ \hfill (Higman-Sims group)

\smallskip

\item[(13)] $Co_3$, $v=276$. \hfill (smallest Conway group)
\end{enumerate}

\medskip

For required basic properties of the listed groups, we refer, e.g.,
to~\cite{Atlas1985,HupI1967,KlLi1990,Suz1962,Ward1966}.

We will now indicate some helpful combinatorial tools on which we
rely in the sequel. Let $r$ (respectively $\lambda_2$) denote the
total number of blocks incident with a given point (respectively
pair of distinct points), and let all further parameters be as
defined at the beginning of Section~\ref{intro}.

Obvious is the subsequent fact.

\begin{lem}\label{divprop}
Let $\D=(X,\B,I)$ be a Steiner \mbox{$t$-design}. If $G \leq
\aut(\D)$ acts flag-transitively on $\D$, then, for any $x \in X$,
the division property
\[r  \bigm|  \left| G_x \right|\] holds.
\end{lem}

Elementary counting arguments give the following standard assertions.

\begin{lem} \label{Comb_t=4}
Let $\D=(X,\B,I)$ be a \mbox{$t$-$(v,k,\lambda)$} design. Then the
following holds:
\begin{enumerate}

\smallskip

\item[(a)] $bk = vr.$

\smallskip

\item[(b)] $\displaystyle{{v \choose t} \lambda = b {k \choose
t}.}$

\smallskip

\item[(c)] $r(k-1)=\lambda_2(v-1)$ for $t \geq 2$, where
$\displaystyle{\lambda_2=\lambda \frac{{v-2 \choose t-2}}{{k-2
\choose t-2}}.}$

\smallskip

\item[(d)] In particular, if $t=4$, then $(k-2)(k-3) \mid (v-2)(v-3).$
\end{enumerate}
\end{lem}

\medskip

For non-trivial Steiner \mbox{$t$-designs} lower bounds for $v$ in
terms of $k$ and $t$ are known.

\begin{prop}{\em (cf.~\cite{Cam1976}).}\label{Cam}
If $\D=(X,\B,I)$ is a non-trivial Steiner \mbox{$t$-design}, then
the following holds:
\begin{enumerate}

\smallskip

\item[(a)] $v\geq (t+1)(k-t+1).$

\smallskip

\item[(b)] $v-t+1 \geq (k-t+2)(k-t+1)$ for $t>2$. If equality
holds, then
\smallskip
$(t,k,v)=(3,4,8),(3,6,22),(3,12,112),(4,7,23)$, or $(5,8,24)$.
\end{enumerate}
\end{prop}

We note that (a) is stronger for $k<2(t-1)$, while (b) is stronger
for $k>2(t-1)$. For $k=2(t-1)$ both assert that $v \geq t^2-1$.

As we are in particular interested in the case when $t=4$, we deduce
from (b) the following upper bound for the positive integer $k$.

\begin{Cor}\label{Cameron_t=4}
Let $\D=(X,\B,I)$ be a non-trivial Steiner \mbox{$4$-design}. Then
\[k \leq \bigl\lfloor \sqrt{v} + \textstyle{\frac{5}{2}} \bigr\rfloor.\]
\end{Cor}

\smallskip

\begin{rem} \label{equa_t=4}
\emph{If \mbox{$G \leq \Aut(\D)$} acts flag-transitively on any
Steiner \mbox{$4$-design} $\D$, then applying
Proposition~\ref{flag2trs} and Lemma~\ref{Comb_t=4}~(b) yields the
equation
\[b=\frac{v(v-1)(v-2)(v-3)}{k(k-1)(k-2)(k-3)}=\frac{v(v-1)
\left|G_{xy}\right|}{\left| G_B \right|},\] where $x$ and $y$ are
two distinct points in $X$ and $B$ is a block in $\B$, and thus
\[(v-2)(v-3) = (k-1)(k-2)(k-3) \frac{\left|G_{xy}
\right|}{\left|G_{xB}\right|} \;\, \mbox{if} \;\, x \in B.\]}
\end{rem}

\medskip

Finally, we assert some lemmas which will be required in
Subsection~\ref{almost simple type} of the proof of the Main
Theorem. Let $q$ be a prime power $p^e$, and $U$ a subgroup of
$PSL(2,q)$. Furthermore, let $N_l$ denote the number of orbits of
length $l$ and let \mbox{$n=(2,q-1)$}.
In~\cite[Ch.\,5]{Hu_Habil2005}, we have in particular determined the
orbit-lengths from the action of subgroups of $PSL(2,q)$ on the
points of the projective line. For the list of subgroups of
$PSL(2,q)$, we thereby refer to~\cite[Ch.\,12,\,p.\,285f.]{Dick1901}
or~\cite[Ch.\,II,\,Thm.\,8.27]{HupI1967}.

\medskip

\begin{lem}\label{PSL_cyc}
Let $U$ be the cyclic group of order $c$ with $c \mid \frac{q \pm
1}{n}$. Then, we have
\begin{enumerate}
\item[(a)] if $c \mid \frac{q+1}{n}$, then $N_c=(q+1)/c$,
\item[(b)] if $c \mid \frac{q-1}{n}$, then $N_1=2$ and $N_c=(q-1)/c$.
\end{enumerate}
\end{lem}

\begin{lem}\label{PSL_dihed}
Let $U$ be the dihedral group of order $2c$ with $c \mid \frac{q \pm
1}{n}$. Then
\begin{enumerate}
\item[(i)] for \mbox{$q \equiv 1$ $($\emph{mod} $4)$}, we have
\begin{enumerate}
\item[(a)] if $c \mid \frac{q+1}{2}$, then $N_c=2$ and $N_{2c}=(q+1-2c)/(2c)$,
\item[(b)] if $c \mid \frac{q-1}{2}$, then $N_2=1$, $N_c=2$, and $N_{2c}=(q-1-2c)/(2c)$,
unless $c=2$, in which case $N_2=3$ and $N_4=(q-5)/4$,
\end{enumerate}
\item[(ii)] for \mbox{$q \equiv 3$ $($\emph{mod} $4)$}, we have
\begin{enumerate}
\item[(a)] if $c \mid \frac{q+1}{2}$, then $N_{2c}=(q+1)/(2c)$,
\item[(b)] if $c \mid \frac{q-1}{2}$, then $N_2=1$ and
$N_{2c}=(q-1)/(2c)$,
\end{enumerate}
\item[(iii)] for \mbox{$q \equiv 0$ $($\emph{mod} $2)$}, we have
\begin{enumerate}
\item[(a)] if $c \mid q+1$, then $N_c=1$ and $N_{2c}=(q+1-c)/(2c)$,
\item[(b)] if $c \mid q-1$, then $N_2=1$, $N_c=1$, and $N_{2c}=(q-1-c)/(2c)$.
\end{enumerate}
\end{enumerate}
\end{lem}

\begin{lem}\label{PSL_elAb}
Let $U$ be the elementary Abelian group of order $\q \mid q$. Then,
we have $N_1=1$ and $N_{\q}=q/ \q$.
\end{lem}

\begin{lem}\label{PSL_semi}
Let $U$ be a semi-direct product of the elementary Abelian group of
order $\q \mid q$ and the cyclic group of order $c$ with $c \mid
\q-1$ and $c \mid q-1$. Then, we have $N_1=1$, $N_{\q}=1$, and $N_{c
\q}=(q- \q)/(c \q)$.
\end{lem}

\begin{lem}\label{PSL_kl.PSL}
Let $U$ be $PSL(2,\q)$ with $\q^m = q$, $m \geq 1$. Then, we have
$N_{\q+1}=1$, $N_{\q (\q-1)}=1$ if $m$ is even, and all other orbits
are regular.
\end{lem}

\begin{lem}\label{PSL_kl.PGL}
Let $U$ be $PGL(2,\q)$ with $\q^m = q$, $m > 1$ even. Then, we have
$N_{\q+1}=1$, $N_{\q(\q-1)}=1$, and all other orbits are regular.
\end{lem}

\begin{lem}\label{PSL_A_4}
Let $U$ be isomorphic to $A_4$. Then
\begin{enumerate}
\item[(i)] for \mbox{$q \equiv 1$ $($\emph{mod} $4)$}, we have
\begin{enumerate}
\item[(a)] if $3 \mid \frac{q+1}{2}$, then $N_{6}=1$ and $N_{12}=(q-5)/12$,
\item[(b)] if $3 \mid \frac{q-1}{2}$, then $N_4=2$, $N_{6}=1$, and $N_{12}=(q-13)/12$,
\item[(c)] if $3 \mid q$, then $N_4=1$, $N_{6}=1$, and $N_{12}=(q-9)/12$,
\end{enumerate}
\item[(ii)] for \mbox{$q \equiv 3$ $($\emph{mod} $4)$}, we have
\begin{enumerate}
\item[(a)] if $3 \mid \frac{q+1}{2}$, then $N_{12}=(q+1)/12$,
\item[(b)] if $3 \mid \frac{q-1}{2}$, then $N_4=2$ and $N_{12}=(q-7)/12$,
\item[(c)] if $3 \mid q$, then $N_4=1$ and $N_{12}=(q-3)/12$,
\end{enumerate}
\item[(iii)] for $q=2^e$, \mbox{$e \equiv 0$ $($\emph{mod} $2)$}, we have
$N_{1}=1$, $N_{4}=1$, and \\$N_{12}=(q-4)/12$.
\end{enumerate}
\end{lem}

\pagebreak

\begin{lem}\label{PSL_S_4}
Let $U$ be isomorphic to $S_4$. Then
\begin{enumerate}
\item[(i)] for \mbox{$q \equiv 1$ $($\emph{mod} $8)$}, we have
\begin{enumerate}
\item[(a)] if $3 \mid \frac{q+1}{2}$, then $N_{6}=1$, $N_{12}=1$, and $N_{24}=(q-17)/24$,
\item[(b)] if $3 \mid \frac{q-1}{2}$, then $N_{6}=1$, $N_{8}=1$, $N_{12}=1$, and \\$N_{24}=(q-25)/24$,
\item[(c)] if $3 \mid q$, then $N_4=1$, $N_{6}=1$, and $N_{24}=(q-9)/24$,
\end{enumerate}
\item[(ii)] for \mbox{$q \equiv -1$ $($\emph{mod} $8)$}, we have
\begin{enumerate}
\item[(a)] if $3 \mid \frac{q+1}{2}$, then $N_{24}=(q+1)/24$,
\item[(b)] if $3 \mid \frac{q-1}{2}$, then $N_8=1$ and
$N_{24}=(q-7)/24$.
\end{enumerate}
\end{enumerate}
\end{lem}

\begin{lem}\label{PSL_A_5}
Let $U$ be isomorphic to $A_5$. Then
\begin{enumerate}
\item[(i)] for \mbox{$q \equiv 1$ $($\emph{mod} $4)$}, we have
\begin{enumerate}
\item[(a)] if $q=5^e$, \mbox{$e \equiv 1$ $($\emph{mod} $2)$}, then $N_{6}=1$ and $N_{60}=(q-5)/60$,
\item[(b)] if $q=5^e$, \mbox{$e \equiv 0$ $($\emph{mod} $2)$}, then $N_{6}=1$, $N_{20}=1$, and \\ $N_{60}=(q-25)/60$,
\item[(c)] if $15 \mid \frac{q+1}{2}$, then $N_{30}=1$ and $N_{60}=(q-29)/60$,
\item[(d)] if $3 \mid \frac{q+1}{2}$ and $5 \mid \frac{q-1}{2}$, then $N_{12}=1$, $N_{30}=1$, and \\ $N_{60}=(q-41)/60$,
\item[(e)] if $3 \mid \frac{q-1}{2}$ and $5 \mid \frac{q+1}{2}$, then $N_{20}=1$, $N_{30}=1$, and \\ $N_{60}=(q-49)/60$,
\item[(f)] if $15 \mid \frac{q-1}{2}$, then $N_{12}=1$, $N_{20}=1$, $N_{30}=1$, and \\ $N_{60}=(q-61)/60$,
\item[(g)] if $3 \mid q$ and $5 \mid \frac{q+1}{2}$, then $N_{10}=1$ and $N_{60}=(q-9)/60$,
\item[(h)] if $3 \mid q$ and $5 \mid \frac{q-1}{2}$, then $N_{10}=1$, $N_{12}=1$, and \\$N_{60}=(q-21)/60$,
\end{enumerate}
\item[(ii)] for \mbox{$q \equiv 3$ $($\emph{mod} $4)$}, we have
\begin{enumerate}
\item[(a)] if $15 \mid \frac{q+1}{2}$, then $N_{60}=(q+1)/60$,
\item[(b)] if $3 \mid \frac{q+1}{2}$ and $5 \mid \frac{q-1}{2}$, then $N_{12}=1$ and $N_{60}=(q-11)/60$,
\item[(c)] if $3 \mid \frac{q-1}{2}$ and $5 \mid \frac{q+1}{2}$, then $N_{20}=1$ and $N_{60}=(q-19)/60$,
\item[(d)] if $15 \mid \frac{q-1}{2}$, then $N_{12}=1$, $N_{20}=1$, and
$N_{60}=(q-31)/60$.
\end{enumerate}
\end{enumerate}
\end{lem}

\bigskip


\section{Proof of the Main Theorem}

Using the notation as before, \textbf{let $\D=(X,\B,I)$ be a
non-trivial Steiner \mbox{$4$-design} with \mbox{$G \leq \Aut(\D)$}
acting flag-transitively on $\D$ throughout the proof}. We recall
that due to Proposition~\ref{flag2trs}, we may restrict ourselves to
the consideration of the finite \mbox{$2$-transitive} permutation
groups listed in Section~\ref{Prelim}. Clearly, in the following we
may assume that $k>4$ as trivial Steiner \mbox{$4$-designs} are
excluded. For each of the Cases (A)(5)-(8) and (B)(8),(11)-(13) we
have only a small number of possibilities for $k$ to check, which
can easily be ruled out by hand using Lemma~\ref{divprop},
Lemma~\ref{Comb_t=4}~(b)-(d), and Corollary~\ref{Cameron_t=4}.


\medskip

\subsection{Groups of Automorphisms of Affine
Type}\label{affine typ} \hfill

\smallskip

We first assume that $G$ is of affine type.

\medskip
\emph{Case} $G \leq A \mathit{\Gamma} L(1,v)$, $v=p^d$.
\medskip

As $G$ is point \mbox{$2$-transitive}, we have $\left| G \right|
=v(v-1)a$ with $a \mid d$. Using Lemma~\ref{divprop}, we obtain
\[(p^d-2)(p^d-3) \bigm| a (k-1)(k-2)(k-3) \bigm| d (k-1)(k-2)(k-3),\]
and hence in particular
\[(p^d-2)(p^d-3) \leq d(k-1)(k-2)(k-3).\]
But, Proposition~\ref{Cam}~(b) yields
\[p^d-3 \geq (k-2)(k-3),\]
and thus
\[p^d-2 \leq d(k-1).\]
With regard to Corollary~\ref{Cameron_t=4}, this leaves only a very
small number of possibilities for $k$ to check, which can easily be
ruled out by hand using Lemma~\ref{Comb_t=4}~(b) and (c). Therefore,
\mbox{$G \leq \Aut(\D)$} cannot act flag-transitively on any
non-trivial Steiner \mbox{$4$-design} $\D$.

\medskip
\emph{Case}  $G_0 \unrhd SL(\frac{d}{a},p^a)$, $d \geq 2a$.
\smallskip

In the following, let $e_i$ denote the $i$-th basis vector of the
vector space $V=V(\frac{d}{a},p^a)$, and
$\text{\footnotesize{$\langle$}} e_i
\text{\footnotesize{$\rangle$}}$ the \mbox{$1$-dimensional} vector
subspace spanned by $e_i$.

First, let $p^a >3$. For $d=2a$, let
$U=U(\text{\footnotesize{$\langle$}} e_1
\text{\footnotesize{$\rangle$}}) \leq G_0$ denote the subgroup of
all transvections with axis $\text{\footnotesize{$\langle$}} e_1
\text{\footnotesize{$\rangle$}}$. Then $U$ consists of all elements
of the form
\[\begin{pmatrix}
  1 & 0 \\
  c & 1
\end{pmatrix},\; c \in GF(p^a) \; \mbox{arbitrary}.\]
Clearly, $U$ fixes as points only the elements of
$\text{\footnotesize{$\langle$}} e_1
\text{\footnotesize{$\rangle$}}$. Thus, $G_0$ has point-orbits of
length at least $p^a$ outside $\text{\footnotesize{$\langle$}} e_1
\text{\footnotesize{$\rangle$}}$. Let \mbox{$S=\{0,e_1,x,y\}$} be a
\mbox{$4$-subset} of distinct points with $x,y \in
\text{\footnotesize{$\langle$}} e_1
\text{\footnotesize{$\rangle$}}$. Obviously, $U$ fixes the unique
block $B \in \B$ which is incident with $S$. If $B$ contains at
least one point outside $\text{\footnotesize{$\langle$}} e_1
\text{\footnotesize{$\rangle$}}$, then we would obtain $k \geq p^a
+4$, which is not possible as $k \leq p^a+2$ in view of
Corollary~\ref{Cameron_t=4}. Thus $B$ is contained completely in
$\text{\footnotesize{$\langle$}} e_1
\text{\footnotesize{$\rangle$}}$, and as $G$ is flag-transitive, we
conclude that each block lies in an affine line. But, by the
definition of Steiner \mbox{$4$-designs}, any four distinct
non-collinear points must also be incident with a
unique block, a contradiction.\\
For $d \geq 3a$, $SL(\frac{d}{a},p^a)_{e_1}$ and hence also
$G_{0,e_1}$ acts point-transitively on $V \setminus
\text{\footnotesize{$\langle$}} e_1
\text{\footnotesize{$\rangle$}}$. Again, let
\mbox{$S=\{0,e_1,x,y\}$} be a \mbox{$4$-subset} of distinct points
with $x,y \in \text{\footnotesize{$\langle$}} e_1
\text{\footnotesize{$\rangle$}}$. If the unique block $B \in \B$
which is incident with $S$ contains some point outside
$\text{\footnotesize{$\langle$}} e_1
\text{\footnotesize{$\rangle$}}$, then it would already contain all
points outside, thus at least $p^d-p^a+4$ many, which obviously
contradicts Corollary~\ref{Cameron_t=4}. We conclude that $B$ lies
completely in $\text{\footnotesize{$\langle$}} e_1
\text{\footnotesize{$\rangle$}}$, and may proceed as above.

Now, let $p^a=2$. Then $v=2^d$. For $d=3$, we have $v=8$ and $k=5$
by Corollary~\ref{Cameron_t=4}, which is not possible in view of
Lemma~\ref{Comb_t=4}~(c). Therefore, we assume that $d>3$. We remark
that clearly any three distinct points are non-collinear in
$AG(d,2)$ and hence define an affine plane. Let $\mathcal{E}=
\text{\footnotesize{$\langle$}} e_1,e_2
\text{\footnotesize{$\rangle$}}$ denote the \mbox{$2$-dimensional}
vector subspace spanned by $e_1$ and $e_2$. Then
$SL(d,2)_\mathcal{E}$ and hence also $G_{0,\mathcal{E}}$ acts
point-transitively on \mbox{$V \setminus \mathcal{E}$}. If the
unique block $B \in \B$ which is incident with the \mbox{$4$-subset}
\mbox{$\{0,e_1,e_2,e_1+e_2\}$} contains some point outside
$\mathcal{E}$, then it would already contain all points of \mbox{$V
\setminus \mathcal{E}$}, and hence \mbox{$k \geq 2^d-4+4=2^d$}, a
contradiction to Corollary~\ref{Cameron_t=4}. Therefore, $B$ can be
identified with $\mathcal{E}$, and the flag-transitivity of $G$
implies that each block must be an affine plane, a contradiction as
$k>4$. Similar arguments hold for $p^a=3$.

\bigskip
\emph{Case} $G_0 \unrhd Sp(\frac{2d}{a},p^a)$, $d \geq 2a$.
\medskip

First, let $p^a \neq 2.$ The permutation group
$PSp(\frac{2d}{a},p^a)$ on the points of the associated projective
space is a rank $3$ group, and the orbits of the one-point
stabilizer are known
(e.g.~\cite[Ch.\,II,\,Thm.\,9.15\,(b)]{HupI1967}). Thus, $G_0 \unrhd
Sp(\frac{2d}{a},p^a)$ has exactly two orbits on \mbox{$V \setminus
\text{\footnotesize{$\langle$}} x \text{\footnotesize{$\rangle$}}$}
$(0 \neq x \in V)$ of length at least
\[\frac{p^a(p^{2d-2a}-1)}{p^a-1}=\sum_{i=1}^{\frac{2d}{a}-2}p^{ia} >
p^d.\] Let $S=\{0,x,y,z\}$ be a \mbox{$4$-subset} with $y,z \in
\text{\footnotesize{$\langle$}} x
\text{\footnotesize{$\rangle$}}$. If the unique block which is
incident with $S$ contains at least one point of \mbox{$V
\setminus \text{\footnotesize{$\langle$}} x
\text{\footnotesize{$\rangle$}}$}, then we would have $k
> p^d+4$, which is impossible since $k \leq p^d+2$ by
Corollary~\ref{Cameron_t=4}. Therefore, we can argue as in the
previous Case.

Now, let $p^a=2$. Then $v=2^{2d}$. For $d=2$ (here $Sp(4,2) \cong
S_6$ as well-known), Corollary~\ref{Cameron_t=4} yields $k \leq 6$.
But, Lemma~\ref{Comb_t=4}~(d) rules out the cases when $k=5$ or $6$.
Thus, let $d>2$. It is easily seen that there are $2^{2d-1}(2^{2d} -
1)$ hyperbolic pairs in the non-degenerate symplectic space
$V=V(2d,2)$, and by Witt's theorem, $Sp(2d,2)$ is transitive on
these hyperbolic pairs. Let $\{x,y\}$ denote a hyperbolic pair, and
$\mathcal{E}=\text{\footnotesize{$\langle$}} x,y
\text{\footnotesize{$\rangle$}}$ the hyperbolic plane spanned by
$\{x,y\}$. As $\mathcal{E}$ is non-degenerate, we have the
orthogonal decomposition
\[V=\mathcal{E}\perp\mathcal{E}^\perp.\]
Obviously, $Sp(2d,2)_{\{x,y\}}$ stabilizes $\mathcal{E}^\perp$ as a
subspace, which implies that \linebreak \mbox{$Sp(2d,2)_{\{x,y\}}
\cong Sp(2d-2,2)$}. As $\Out(Sp(2d,2))=1$, we have therefore
\[Sp(2d-2,2) \cong Sp(2d,2)_{\{x,y\}} \unlhd Sp(2d,2)_\mathcal{E}=G_{0,\mathcal{E}}.\]
As $Sp(2d-2,2)$ acts transitively on the non-zero vectors of the
\mbox{$(2d-2)$}-dimensional symplectic subspace, the smallest orbit
on \mbox{$V \setminus \mathcal{E}$} under $G_{0,\mathcal{E}}$ has
length at least $2^{2d-2}-1$. If the unique block $B \in \B$ which
is incident with the \mbox{$4$-subset} \mbox{$\{0,x,y,x+y\}$}
contains some point in \mbox{$V \setminus \mathcal{E}$}, then we
would have $k \geq 2^{2d-2}+3$, a contradiction to
Corollary~\ref{Cameron_t=4}. Thus $B$ can be identified with
$\mathcal{E}$, and by the flag-transitivity of $G$, each block must
be an affine plane, yielding a contradiction.

\bigskip
\emph{Case} $G_0 \unrhd G_2(2^a)'$, $d=6a$.
\medskip

First, let $a=1$. Then $v=2^6=64$ and  $k \leq 10$ by
Corollary~\ref{Cameron_t=4}. On the other hand, we have $\left|
G_2(2)' \right|=2^5 \cdot 3^3 \cdot 7$ and $\left| \Out(G_2(2)')
\right| =2$. In view of Lemma~\ref{divprop} this gives
\[r=\frac{63 \cdot 62 \cdot 61}{(k-1)(k-2)(k-3)}
 \Bigm|  \left| G_0 \right|
\Bigm|  2^6 \cdot 3^3 \cdot 7,\] which implies that $k$ is at least
$63$, a contradiction.

Now, let $a>1$. As here $G_2(2^a)$ is simple non-Abelian, it is
sufficient to consider $G_0 \unrhd G_2(2^a)$. The permutation group
$G_2(2^a)$ is of rank $4$, and for $0 \neq x \in V$ the one-point
stabilizer $G_2(2^a)_x$ has exactly three orbits $\mathcal{O}_i$
$(i=1,2,3)$ on \mbox{$V \setminus \text{\footnotesize{$\langle$}} x
\text{\footnotesize{$\rangle$}}$} of length
$2^{3a}-2^a,2^{5a}-2^{3a},2^{6a}-2^{5a}$ (cf., e.g.,~\cite{Asch1987}
or~\cite[Thm.\,3.1]{CaKa1979}). Thus, $G_0$ has exactly three orbits
on \mbox{$V \setminus \text{\footnotesize{$\langle$}} x
\text{\footnotesize{$\rangle$}}$} of length at least $\left|
\mathcal{O}_i \right|.$ Let $S=\{0,x,y,z\}$ be a \mbox{$4$-subset}
with $y,z \in \text{\footnotesize{$\langle$}} x
\text{\footnotesize{$\rangle$}}$. If the unique block $B \in \B$
which is incident with $S$  contains at least one point of \mbox{$V
\setminus \text{\footnotesize{$\langle$}} x
\text{\footnotesize{$\rangle$}}$} in $\mathcal{O}_2$ or
$\mathcal{O}_3$, then we would obtain as above a contradiction to
Corollary~\ref{Cameron_t=4}. Thus, we only have to consider the case
when $B$ contains points of \mbox{$V \setminus
\text{\footnotesize{$\langle$}} x \text{\footnotesize{$\rangle$}}$}
which all lie in $\mathcal{O}_1$. By~\cite{Asch1987}, the orbit
$\mathcal{O}_1$ is exactly known, and we have
\[\mathcal{O}_1 = x\Delta \setminus
\text{\footnotesize{$\langle$}} x \text{\footnotesize{$\rangle$}},\]
where $x\Delta = \{y \in V \mid f(x,y,z)=0 \;\,\mbox{for all} \;\, z
\in V\}$ with an alternating trilinear form $f$ on $V$. Then $B$
consists, apart from elements of $\text{\footnotesize{$\langle$}} x
\text{\footnotesize{$\rangle$}}$, exactly of $\mathcal{O}_1$. Since
$\left| \mathcal{O}_1 \right| \neq 1$, we can choose
$\text{\footnotesize{$\langle$}} \overline{x}
\text{\footnotesize{$\rangle$}} \in x\Delta$ with
$\text{\footnotesize{$\langle$}} \overline{x}
\text{\footnotesize{$\rangle$}} \neq \text{\footnotesize{$\langle$}}
x \text{\footnotesize{$\rangle$}}$. But then, for symmetric reasons,
the \mbox{$4$-subset} $\{0,\overline{x},\overline{y},\overline{z}\}$
with $\overline{y},\overline{z} \in \text{\footnotesize{$\langle$}}
\overline{x} \text{\footnotesize{$\rangle$}}$ must also be incident
with the unique block $B$, a contradiction to the fact that
$\overline{x} \Delta \neq x \Delta$ for
$\text{\footnotesize{$\langle$}} \overline{x}
\text{\footnotesize{$\rangle$}} \neq \text{\footnotesize{$\langle$}}
x \text{\footnotesize{$\rangle$}}$. Consequently, $B$ is contained
completely in $\text{\footnotesize{$\langle$}} x
\text{\footnotesize{$\rangle$}}$, and we may argue as in the Cases
above.

\medskip

\subsection{\mbox{Groups of Automorphisms of Almost Simple
Type}}\label{almost simple type} \hfill

\medskip

We consider now successively those cases where $G$ is of almost
simple type.

\medskip
\emph{Case} $N=A_v$, $v \geq 5$.
We may assume that $v \geq 6$, and then $A_v$ and hence also $G$ is
\mbox{$4$-transitive} and does not act on any non-trivial Steiner
\mbox{$4$-design} $\D$ in view of~\cite[Thm.\,3]{Kant1985}.

\bigskip
\emph{Case} $N=PSL(d,q)$, $d \geq 2$, $v=\frac{q^d-1}{q-1}$, where
$(d,q) \not= (2,2),(2,3)$.
\medskip

We distinguish two subcases:

\medskip
\indent \emph{Subcase} $N=PSL(2,q)$, $v=q+1$, $q=p^e >3$.
\smallskip

Here $\Aut(N)= P \mathit{\Gamma} L (2,q)$, and $\left| G \right| =
(q+1)q \frac{(q-1)}{n}a$ with $n=(2,q-1)$ and $a \mid ne$. We may
assume that $q \geq 5$.

\textbf{First, we assume that $N=G$.} Then, by
Remark~\ref{equa_t=4}, we obtain
\begin{equation}\label{Eq-0}
(q-2) \left| PSL (2,q)_{0B} \right| \cdot n  = (k-1)(k-2)(k-3)
\end{equation}
which is equivalent to
\begin{equation}\label{Eq-0-equiv}
(q-2) \left| PSL (2,q)_{0B} \right| \cdot n +6=k(k^2-6k+11).
\end{equation}
Thus, we have in particular
\begin{equation}\label{Eq-1}
k \bigm| (q-2) \left| PSL (2,q)_{0B} \right| \cdot n +6.
\end{equation}
Since $PSL(2,q)_B$ acts transitively on the points of $B$, we have
\begin{equation}\label{Eq-k}
k=\left|0^{PSL(2,q)_B}\right| = \big|PSL(2,q)_B : PSL(2,q)_{0B}
\big|.
\end{equation}
Let us first consider the case~(1.1) when $\left|PSL(2,q)_{0B}
\right|=1$. a) If $q$ is even, then $k \mid q+4$ by
property~(\ref{Eq-1}). On the other hand, using
equation~(\ref{Eq-k}), we have $k=\left| PSL (2,q)_B \right| \bigm|
\left| PSL (2,q) \right|=q^3-q$. As
\[(q^3-q,q+4)=(60,q+4)=4\cdot
(15,2^{e-2}+1)=\left\{\begin{array}{lll}
    4,\;\, & \mbox{if} \;\, e \;\,\mbox{is even and}\;\,4 \nmid e\\
    4 \cdot 3,\;\,& \mbox{if} \;\,e\;\,\mbox{is odd}\\
    4 \cdot 5,\;\, & \mbox{if} \;\,4 \mid e
\end{array} \right.\] the possible values for $k$ can immediately be
ruled out by hand using equation~(\ref{Eq-0}).

b) If $q$ is odd, we have $k=\left| PSL (2,q)_B \right| \bigm|
2(q+1)$ due to property~(\ref{Eq-1}) and equation~(\ref{Eq-k}).
Examining the list of subgroups of $PSL(2,q)$
(cf.~\cite[Ch.\,12,\,p.\,285f.]{Dick1901}
or~\cite[Ch.\,II,\,Thm.\,8.27]{HupI1967}), we have to consider the
following possibilities:

\begin{enumerate}
\item[(1.1.b)](i) $PSL(2,q)_B$ is conjugate to a cyclic
subgroup of order $c$ with $c \mid \frac{q + 1}{2}$ of $PSL(2,q)$,
and $k=c$.

\item[(1.1.b)](ii) $PSL(2,q)_B$ is conjugate to a dihedral subgroup of order
$2c$ with $c \mid \frac{q + 1}{2}$ of $PSL(2,q)$, and $k=2c$.

\item[(1.1.b)](iii) $PSL(2,q)_B$ is conjugate to $A_4$, and $k=12$.

\item[(1.1.b)](iv) $PSL(2,q)_B$ is conjugate to $S_4$, and $k=24$.

\item[(1.1.b)](v) $PSL(2,q)_B$ is conjugate to $A_5$, and $k=60$.

\end{enumerate}

(1.1.b)(i): By equation~(\ref{Eq-0}), we have
\[c \bigm| \frac{q+1}{2}=\frac{(c-1)(c-2)(c-3)+6}{4}=\frac{c(c^2-6c+11)}{4}.\]
Since $4$ does not divide $c^2-6c+11$, this is impossible.

(1.1.b)(ii): Using equation~(\ref{Eq-0}), we obtain
\[c \bigm| \frac{q+1}{2}=\frac{(2c-1)(c-1)(2c-3)+3}{2},\]
which is not possible as $(2c-1)(c-1)(2c-3)+3=4c^3-12c^2+11c \equiv
c$ (mod $2c$).

(1.1.b)(iii)-(v): For each given value of $k$, equation~(\ref{Eq-0})
gives in each subcase that $q$ is not a prime power.

We consider now the case~(1.2) when $\left| PSL (2,q)_{0B}
\right|=2.$ a) If $q$ is even, then we have $k=\frac{\left|
PSL(2,q)_B \right|}{2} \bigm| 2(q+1)$ due to property~(\ref{Eq-1})
and equation~(\ref{Eq-k}). Considering the list of subgroups of
$PSL(2,q)$, we have the following possibilities:

\begin{enumerate}
\item[(1.2.a)](i) $PSL(2,q)_B$ is conjugate to a cyclic
subgroup of order $c$ with $c \mid q + 1$ of $PSL(2,q)$, and
$k=\frac{c}{2}$.

\item[(1.2.a)](ii) $PSL(2,q)_B$ is conjugate to a dihedral subgroup of order
$2c$ with $c \mid q + 1$ of $PSL(2,q)$, and $k=c$.

\item[(1.2.a)](iii) $PSL(2,q)_B$ is conjugate to $PSL(2,\q)$ with
$\q \mid 4$, and $k=30$.

\item[(1.2.a)](iv) $PSL(2,q)_B$ is conjugate to $A_4$, and $k=6$.
\end{enumerate}

(1.2.a)(i),(iii),(iv): In view of
Lemmas~\ref{PSL_cyc},~\ref{PSL_kl.PGL},
respectively~\ref{PSL_A_4}~(iii), clearly $k$ cannot take the given
values.

(1.2.a)(ii): Considering equation~(\ref{Eq-0}), for $k=c>4$ clearly
the right hand side of the equation is divisible by $8$, but not the
left hand side.

b) If $q$ is odd, then $k \mid 2(2q-1)$ by property~(\ref{Eq-1}). On
the other hand, equation~(\ref{Eq-k}) yields $k=\frac{\left| PSL
(2,q)_B \right|}{2} \Bigm| \frac{\left| PSL (2,q)
\right|}{2}=\frac{q^3-q}{4}$. Since
$\big(\frac{q^3-q}{4},2(2q-1)\big) \linebreak = 2 \cdot
\big(\frac{q^3-q}{8},2q-1 \big)=2\cdot (3,q+1)$ only $k=6$ can
occur, and equation~(\ref{Eq-0}) gives then $q=17$. However, it is
known that there does not exist any \mbox{$4$-$(18,6,1)$ design}
(cf.~\cite[Thm.\,6]{Witt1938b}).

Finally, let us consider the case~(1.3) when $\left| PSL (2,q)_{0B}
\right|>2$. Examining the list of subgroups of $PSL(2,q)$ with their
orbits on the projective line (Lemmas~\ref{PSL_cyc}-\ref{PSL_A_5}),
we have to consider the following subcases:

\begin{enumerate}

\item[(1.3)](i) $PSL(2,q)_B$ is conjugate to $S_4$, and $k=6$ or $8$.

\item[(1.3)](ii) $PSL(2,q)_B$ is conjugate to $A_5$, and $k=6,10,12$ or $20$.

\item[(1.3)](iii) $PSL(2,q)_B$ is conjugate to a semi-direct product of an elementary Abelian subgroup of order
$\q \mid q$ with a cyclic subgroup of order $c$ of $PSL(2,q)$ with
$c \mid \q-1$ and $c \mid q-1$, and $k=\q$.

\item[(1.3)](iv) $PSL(2,q)_B$ is conjugate to $PSL(2,\q)$ with $\q^m = q$, $m \geq 1$, and $k=
\q +1$ or $\q(\q-1)$ if $m$ is even.

\item[(1.3)](v) $PSL(2,q)_B$ is conjugate to $PGL(2,\q)$ with $\q^m = q$, $m > 1$ even, and $k=
\q +1$ or $\q(\q-1)$.
\end{enumerate}

(1.3)(i): We may assume that $q$ is odd. Applying
equations~(\ref{Eq-0}) and~(\ref{Eq-k}) yields for $k=6$ that $q$ is
not a prime power, and for $k=8$ that $q=37$, in which case $q\equiv
\pm 1$ (mod $8$) (cf. Lemma~\ref{PSL_S_4}) does not hold.

(1.3)(ii): Again, we may assume that $q$ is odd and consider
equations~(\ref{Eq-0}) and~(\ref{Eq-k}) for the given values of $k$.
We obtain for $k=6$ that $q=5$, which is clearly impossible due to
Corollary~\ref{Cameron_t=4}, for $k=10$ that $q$ is not a prime
power, and for $k=20$ that $q=971$, in which case
Lemma~\ref{Comb_t=4}~(c) gives a contradiction. If $k=12$, then we
get $q=101$. Since $\left| PSL (2,q)_{0B} \right|=5$ by
equation~(\ref{Eq-k}) and $5 \bigm| \frac{q-1}{2}$, $PSL(2,q)_{0B}$
has two distinct fixed points. If one fixed point lies outside $B$,
then clearly $q \equiv 1$ (mod $5$) and hence $k=12$ is not
possible. Thus, we may assume that both fixed points are incident
with $B$. But then, as every non-identity element of $PSL(2,q)$
fixes at most two distinct points, $PSL(2,q)_{0B}$ must fix some
\mbox{$2$-subset} by the definition of Steiner \mbox{$4$-designs},
and hence contains an involution, a contradiction.

(1.3)(iii): We have $\big((q-2) \left| PSL (2,q)_{0B} \right| \cdot
n +6,q \big)=\big(2 \cdot \left| PSL (2,q)_{0B} \right| \cdot n -6,q
\big)$, and property~(\ref{Eq-1}) gives in particular
\begin{equation}\label{Eq_>2_a}
k \bigm| 2 \cdot \left| PSL (2,q)_{0B} \right| \cdot n -6.
\end{equation}
On the other hand, it follows from equation~(\ref{Eq-k}) that
$\left| PSL (2,q)_{0B} \right| \bigm| k-1$. Therefore, we have in
particular
\begin{equation}\label{Eq_>2}
\frac{k-1}{2n}<\frac{k+6}{2n} \leq \left| PSL (2,q)_{0B} \right|
\bigm| k-1.
\end{equation}
a) If $q$ is even, then we deduce that $\left| PSL (2,q)_{0B}
\right| = k-1$. Property~(\ref{Eq_>2_a}) yields that $k \bigm| 2
k-8$, and, as clearly $(2k-8,k)=(8,k)$, only $k=8$ is possible.
Thus, we have $q=32$ in view of equation~(\ref{Eq-0}), which is
impossible by Lemma~\ref{Comb_t=4}~(c).

b) If $q$ is odd, then by property~(\ref{Eq_>2}), we have to
consider the possibilities when $\left| PSL (2,q)_{0B} \right| =
\frac{k-1}{\overline{n}}$ with $\overline{n}=1,2,3$. If $\left| PSL
(2,q)_{0B} \right| = k-1$, then we obtain $k \bigm| 4k-10$ by
property~(\ref{Eq-1}). Clearly, $(4k-10,k)=(10,k)$, but as $k \mid
q$, only $k=5$ is possible. Then, equation~(\ref{Eq-0}) gives $q=5$,
which leads to a contradiction in view of
Corollary~\ref{Cameron_t=4}. For $\left| PSL (2,q)_{0B} \right| =
\frac{k-1}{\overline{2}}$, we have $k \bigm| 2 k-8$ and thus $k=8$
as above, which is impossible as $k \nmid q$. If $\left| PSL
(2,q)_{0B} \right| = \frac{k-1}{\overline{3}}$, then
property~(\ref{Eq-1}) yields $k \bigm| \frac{4k-22}{3}$. Since
$(4k-22,3k)=(22,3k)$ is not divisible by $3$, this is not possible.

(1.3)(iv): We have $\left| PSL (2,q)_{0B} \right| = \frac{\q
(\q-1)}{n}$ if $k=\q+1$. Thus, equation~(\ref{Eq-0}) yields for
$k=\q +1$ that $q=\q$ must hold, which is impossible due to
Corollary~\ref{Cameron_t=4}. For $m>1$ even and $k=\q(\q-1)$, it
follows that $\left| PSL (2,q)_{0B} \right| = \frac{\q+1}{n}$.
Hence, property~(\ref{Eq-1}) gives
\[\q(\q-1) \bigm| (q-2) (\q+1) +6= \q^{m+1}+\q^m-2\q +4.\]
Since $(\q^{m+1}+\q^m-2\q +4,\q)=(4,\q)$ and $k>4$, only the case
when $\q=4$ has to be considered. Thus, $k=12$ and applying
equation~(\ref{Eq-0-equiv}) immediately gives a contradiction.

(1.3)(v): Clearly $n$ does not appear in equations~(\ref{Eq-0})
and~(\ref{Eq-0-equiv}) as well as in property~(\ref{Eq-1}), and we
may argue mutatis mutandis as in subcase~(iv).

\textbf{Now, let us assume that $N<G \leq \Aut(N)$.} We recall that
$q=p^e >3$, and will distinguish in the following the cases $p>2$
and $p=2$.

\textbf{First, let $p>2$.} We define $G^*=G \cap (PSL(2,q) \rtimes
\text{\footnotesize{$\langle$}} \tau_\alpha
\text{\footnotesize{$\rangle$}})$ with $\tau_\alpha \in$
Sym$(GF(p^e) \cup \{\infty\}) \cong S_v$ of order $e$ induced by the
Frobenius automorphism $\alpha : GF(p^e) \longrightarrow GF(p^e),\,
x \mapsto x^p$. Then, by Dedekind's law, we can write
\begin{equation}\label{Eq_G^*0}
G^*= PSL(2,q) \rtimes (G^* \cap \text{\footnotesize{$\langle$}}
\tau_\alpha \text{\footnotesize{$\rangle$}}).
\end{equation}
Defining $P \Si L(2,q)= PSL(2,q) \rtimes
\text{\footnotesize{$\langle$}} \tau_\alpha
\text{\footnotesize{$\rangle$}}$, it can easily by calculated that
$P \Si L(2,q)_{0,1,\infty} = \text{\footnotesize{$\langle$}}
\tau_\alpha \text{\footnotesize{$\rangle$}}$, and
$\text{\footnotesize{$\langle$}} \tau_\alpha
\text{\footnotesize{$\rangle$}}$ has precisely $p+1$ distinct fixed
points (cf., e.g.,~\cite[Ch.\,6.4,\,Lemma\,2]{Demb1968}). As $p>2$,
we conclude therefore that $G^* \cap \text{\footnotesize{$\langle$}}
\tau_\alpha \text{\footnotesize{$\rangle$}} \leq G^*_{0B}$ for some
appropriate, unique block $B \in \B$ by the definition of Steiner
\mbox{$4$-designs}. Furthermore, clearly $PSL(2,q) \cap (G^* \cap
\text{\footnotesize{$\langle$}} \tau_\alpha
\text{\footnotesize{$\rangle$}}) =1.$ Thus, if we assume that
\mbox{$G^* \leq \Aut(\D)$} acts already flag-transitively on $\D$,
then we obtain $\left| (0,B)^{G^*} \right|=\left|
(0,B)^{PSL(2,q)}\right|=bk$ in view of Remark~\ref{equa_t=4}. Hence,
$PSL(2,q)$ must also act flag-transitively on $\D$, which has
already been considered. Therefore, let us assume that \mbox{$G^*
\leq \Aut(\D)$} does not act flag-transitively on $\D$. Then $\big|G
: G^* \big|=2$ and $G^*$ has exactly two orbits of equal length on
the set of flags. It follows for the orbit containing the flag
$(0,B)$ that $\left| (0,B)^{G^*} \right|=\left|
(0,B)^{PSL(2,q)}\right|=\frac{bk}{2}$. As $PSL(2,q)$ is normal in
$G$, we have under $PSL(2,q)$ also precisely one further orbit of
equal length on the set of flags. Then, proceeding similarly to the
case $N=G$ for each orbit on the set of flags, we obtain
(representative for the orbit containing the flag $(0,B)$) that
\begin{equation}\label{Eq-0_N<G}
\frac{(q-2) \left| PSL (2,q)_{0B} \right| \cdot n}{2}  =
(k-1)(k-2)(k-3)
\end{equation}
which is equivalent to
\begin{equation}\label{Eq-0_N<G-equiv}
\frac{(q-2) \left| PSL (2,q)_{0B} \right| \cdot n}{2}
+6=k(k^2-6k+11).
\end{equation}
Hence, we have in particular
\begin{equation}\label{Eq-1_N<G}
k \Bigm| \frac{(q-2) \left| PSL (2,q)_{0B} \right| \cdot n}{2} +6.
\end{equation}
Since $PSL(2,q)_B$ can have one or two orbits of equal length on the points of $B$, we have
\begin{equation}\label{Eq-k_N<G}
k\;\,\mbox{or}\;\,\frac{k}{2}=\left|0^{PSL(2,q)_B}\right| =
\big|PSL(2,q)_B : PSL(2,q)_{0B} \big|.
\end{equation}
Let us recall that here $q$ is always odd. First considering the
case~(2.1) when $\left|PSL(2,q)_{0B} \right|=1$ yields immediately a
contradiction to equation~(\ref{Eq-0_N<G}). Let us now observe the
case~(2.2) when $\left| PSL (2,q)_{0B} \right|=2$. We have $k \bigm|
2(q+1)$ in view of property~(\ref{Eq-1_N<G}), and $k=\left|
PSL(2,q)_B \right|$ or $\frac{\left| PSL(2,q)_B \right|}{2}$ by
equation~(\ref{Eq-k_N<G}). a) For $k=\left| PSL(2,q)_B \right|$,
clearly equation~(\ref{Eq-0_N<G}) with $\left| PSL (2,q)_{0B}
\right|=2$ is equivalent to equation~(\ref{Eq-0}) with $\left| PSL
(2,q)_{0B} \right|=1$, and thus we can argue exactly as in
case~(1.1.b). b) For $k=\frac{\left| PSL(2,q)_B \right|}{2}$, we
have to consider the following subgroups of $PSL(2,q)$:

\begin{enumerate}
\item[(2.2.b)](i) $PSL(2,q)_B$ is conjugate to a cyclic
subgroup of order $c$ with $c \mid \frac{q + 1}{2}$ of $PSL(2,q)$,
and $k=\frac{c}{2}$.

\item[(2.2.b)](ii) $PSL(2,q)_B$ is conjugate to a dihedral subgroup of order
$2c$ with $c \mid \frac{q + 1}{2}$ of $PSL(2,q)$, and $k=c$.

\item[(2.2.b)](iii) $PSL(2,q)_B$ is conjugate to $A_4$, and $k=6$.

\item[(2.2.b)](iv) $PSL(2,q)_B$ is conjugate to $S_4$, and $k=12$.

\item[(2.2.b)](v) $PSL(2,q)_B$ is conjugate to $A_5$, and $k=30$.
\end{enumerate}

(2.2.b)(i): Obviously, $k$ cannot take the given value due to
Lemma~\ref{PSL_cyc}.

(2.2.b)(ii): It follows from equation~(\ref{Eq-0_N<G}) that
\[c \bigm| \frac{q+1}{2}=\frac{(c-1)(c-2)(c-3)+6}{4}=\frac{c(c^2-6c+11)}{4}.\]
As $4$ does not divide $c^2-6c+11$, this is impossible.

(2.2.b)(iii)-(v): In view of equation~(\ref{Eq-0_N<G}), we obtain in
subcase~(iii) that $q=32$ which is not possible, and in each of the
other subcases that $q$ is not a prime power.

We examine finally the case~(2.3) when $\left| PSL (2,q)_{0B}
\right|>2.$ Combining equations~(\ref{Eq-0_N<G})
and~(\ref{Eq-k_N<G}), we obtain
\begin{equation}\label{Eq-2orbits_1}
\frac{(q-2) \left| PSL(2,q)_B \right| \cdot n}{2}= k (k-1) (k-2) (k-3)
\end{equation}
\[\mbox{with} \;\, k=\left| 0^{PSL(2,q)_B} \right| = \frac{\left| PSL(2,q)_B \right|}{\left| PSL(2,q)_{0B}\right|},
\;\, \mbox{or}\]
\begin{equation}\label{Eq-2orbits_2}
(q-2) \left| PSL(2,q)_B \right| \cdot n = k (k-1) (k-2) (k-3)
\end{equation}
\[\mbox{with} \;\, k= 2 \cdot \left| 0^{PSL(2,q)_B} \right| = 2 \cdot \frac{\left| PSL(2,q)_B \right|}{\left| PSL(2,q)_{0B}\right|}.\]
In view of the subgroups of $PSL(2,q)$ with their orbits on the
projective line (Lemmas~\ref{PSL_cyc}-\ref{PSL_A_5}), we have the
following possibilities:

\begin{enumerate}
\item[(2.3)](i) $PSL(2,q)_B$ is conjugate to $A_4$, and $k= 2 \cdot \left| 0^{PSL(2,q)_B} \right|=8$.

\item[(2.3)](ii) $PSL(2,q)_B$ is conjugate to $S_4$, and $k=\left| 0^{PSL(2,q)_B} \right|= 6$ or
$8$,\\respectively $k= 2 \cdot \left| 0^{PSL(2,q)_B} \right|=8,12$
or $16$.

\item[(2.3)](iii) $PSL(2,q)_B$ is conjugate to $A_5$, and $k=\left| 0^{PSL(2,q)_B} \right|=6,10,12$ or $20,$
respectively $k= 2 \cdot \left| 0^{PSL(2,q)_B} \right|=12,20,24$ or
$40$.

\item[(2.3)](iv) $PSL(2,q)_B$ is conjugate to a semi-direct product of an elementary Abelian subgroup of order
$\q \mid q$ with a cyclic subgroup of order $c$ of $PSL(2,q)$ with
$c \mid \q-1$ and $c \mid q-1$, and $k=\left| 0^{PSL(2,q)_B}
\right|=\q$, respectively $k= 2 \cdot \left| 0^{PSL(2,q)_B}
\right|=2 \q$.

\item[(2.3)](v) $PSL(2,q)_B$ is conjugate to $PSL(2,\q)$ with $\q^m = q$, $m \geq 1$,
and $k=\left| 0^{PSL(2,q)_B} \right|=\q+1$, or $\q(\q-1)$ if $m$ is
even, respectively $k= 2 \cdot \left| 0^{PSL(2,q)_B} \right|=2
(\q+1)$, or $2 \q(\q-1)$ if $m$ is even.

\item[(2.3)](vi) $PSL(2,q)_B$ is conjugate to $PGL(2,\q)$ with $\q^m = q$, $m > 1$ even, and
$k=\left| 0^{PSL(2,q)_B} \right|=\q+1$ or $\q(\q-1)$, respectively
$k= 2 \cdot \left| 0^{PSL(2,q)_B} \right|=2 (\q+1)$ or $2 \q(\q-1)$.
\end{enumerate}

(2.3)(i): By equation~(\ref{Eq-2orbits_2}), we obtain that $q$ is
not a prime power.

(2.3)(ii): First, applying equation~(\ref{Eq-2orbits_1}) yields for
$k=6$ that $q=17$, which can be excluded since there does not exist
any \mbox{$4$-$(18,6,1)$ design} as already mentioned, and for $k=8$
that $q$ is not a prime power. Using equation~(\ref{Eq-2orbits_2})
gives for $k=8$ that $q=37$, in which case $q\equiv \pm 1$ (mod $8$)
(cf. Lemma~\ref{PSL_S_4}) does not hold, and for $k=12$ and $16$
that $q$ is not a prime power in each case.

(2.3)(iii): Observing first equation~(\ref{Eq-2orbits_1}) yields for
each given value of $k$ that $q$ would be even. Now, applying
equation~(\ref{Eq-2orbits_2}) gives for $k=12$ the prime $q=101$,
which is impossible since according to Lemma~\ref{PSL_A_5} we only
have orbits of length $6$ when $p=5$, and for $k=20$ that $q=971$,
in which case Lemma~\ref{Comb_t=4}~(c) gives a contradiction. For
$k=24$ and $40$, we obtain in each case that $q$ is not a prime
power.

(2.3)(iv): Let $k=\q$. As $\big(\frac{(q-2) \left| PSL (2,q)_{0B}
\right| \cdot n}{2} +6,q \big)=\big(\left| PSL (2,q)_{0B} \right|
\cdot n -6,q \big)$, property~(\ref{Eq-1_N<G}) implies that
\begin{equation}\label{Eq_>2_2orbits}
k \bigm| \left| PSL (2,q)_{0B} \right| \cdot n -6.
\end{equation}
On the other hand, as $k=\left|0^{PSL(2,q)_B}\right| =
\big|PSL(2,q)_B : PSL(2,q)_{0B} \big|$, we have $\left| PSL
(2,q)_{0B} \right|=c \bigm| k-1$. Thus,  in particular
\[\frac{k-1}{n}<\frac{k+6}{n} \leq \left| PSL (2,q)_{0B} \right|
\bigm| k-1,\] and hence $\left| PSL (2,q)_{0B} \right| = k-1$ as $q$
is odd. But, property~(\ref{Eq_>2_2orbits}) gives $k \bigm| 2 k-8$,
and as clearly $(2k-8,k)=(8,k)$, it would follow that $k=8$, which
is impossible since $q$ is odd. For $k=2\q$, it follows from
equation~(\ref{Eq-2orbits_2}) that
\[(q-2) n =  4 \cdot \frac{(\q-1)}{c} (2\q-1) (2\q-3),\]
which gives a contradiction as clearly the left hand side of the
equation is not divisible by $4$.

(2.3)(v): For $k= \q+1$, it follows from
equation~(\ref{Eq-2orbits_1}) that $q=2(\q-1)$, which is obviously
impossible for $\q >2$. If $m>1$ even and $k=\q(\q-1)$, then we have
\[(q-2)(\q+1)=2(\q^2-\q-1)(\q^2-\q-2)(\q^2-\q-3)\]
in view of equation~(\ref{Eq-2orbits_1}), which is impossible since
$\q^2-\q-2=(\q+1)(\q-2)$ and $q$ is odd. If $k=2(\q+1)$, then
applying equation~(\ref{Eq-2orbits_2}) gives
\[(q-2)(\q-1)=4(2\q+1)(2\q-1).\]
Clearly $(4(2\q+1)(2\q-1),\q-1)=(\q-1,12)$, and the few
possibilities for $\q$ can be ruled as $q=\q^m$. For $m>1$ even and
$k=2\q(\q-1)$, equation~(\ref{Eq-2orbits_2}) yields
\[(q-2)(\q+1)=4(2\q^2-2\q-1)(\q^2-\q-1)(2\q^2-2\q-3).\]
Here $(4(2\q^2-2\q-1)(\q^2-\q-1)(2\q^2-2\q-3),\q+1)=(\q+1,12)$, and
the few possibilities for $\q$ can be ruled out again.

(2.3)(vi): Clearly $n$ does not appear in
equations~(\ref{Eq-2orbits_1}) and~(\ref{Eq-2orbits_2}), and we may
argue mutatis mutandis as in subcase~(v).

\textbf{Now, let $p=2$.} Then clearly $N=PSL(2,q)=PGL(2,q)$, and
$\Aut(N)=P\Si L(2,q)$. If we assume that
$\text{\footnotesize{$\langle$}} \tau_\alpha
\text{\footnotesize{$\rangle$}} \leq P \Si L (2,q)_{0B}$ for some
appropriate, unique block $B \in \B$, then $PSL(2,q)$ must also be
flag-transitive, which has already been considered. Hence, we may
assume that $\text{\footnotesize{$\langle$}} \tau_\alpha
\text{\footnotesize{$\rangle$}} \nleq P\Si L(2,q)_{0B}$. Let $s$ be
a prime divisor of $e=\left| \text{\footnotesize{$\langle$}}
\tau_\alpha \text{\footnotesize{$\rangle$}} \right|$. As the normal
subgroup $H:=(P \Si L (2,q)_{0,1,\infty})^s \leq
\text{\footnotesize{$\langle$}} \tau_\alpha
\text{\footnotesize{$\rangle$}}$ of index $s$ fixes at least four
distinct points, we have $G \cap H \leq G_{0B}$ for some
appropriate, unique block $B \in \B$ by the definition of Steiner
\mbox{$4$-designs}. It can then be deduced that $e=s^u$ for some $u
\in \Z^+$, since if we assume for $G= P \Si L(2,q)$ that there
exists a further prime divisor $\overline{s}$ of $e$ with
$\overline{s} \neq s$ then $\overline{H}:=(P \Si L
(2,q)_{0,1,\infty})^{\overline{s}} \leq
\text{\footnotesize{$\langle$}} \tau_\alpha
\text{\footnotesize{$\rangle$}}$ and $H$ are both subgroups of $P\Si
L(2,q)_{0B}$ by the flag-transitivity of $P \Si L (2,q)$, and hence
$\text{\footnotesize{$\langle$}} \tau_\alpha
\text{\footnotesize{$\rangle$}} \leq P\Si L(2,q)_{0B}$, a
contradiction. Furthermore, as $\text{\footnotesize{$\langle$}}
\tau_\alpha \text{\footnotesize{$\rangle$}}\nleq P \Si L
(2,q)_{0B}$, we may, by applying Dedekind's law, assume that
\[G_{0B} = PSL(2,q) _{0B} \rtimes (G \cap H).\]
Clearly, $k=\left|0^{G_B}\right| = \big|G_B : G_{0B} \big|$. If $G=
PSL(2,q) \rtimes (G \cap H)$, then $PSL(2,q)$ itself must be
flag-transitive. Therefore, we may assume that $G = P \Si L (2,q)$.
Thus, by Remark~\ref{equa_t=4}, we obtain
\begin{equation}\label{condB}
(q-2) \left| PSL (2,q)_{0B} \right| =(k-1)(k-2)(k-3) s
\end{equation}
\[\mbox{with} \;\, \left| PSL(2,q)_{0B}\right| =\frac{\left| PSL(2,q)_B
\right|}{k} \cdot \left\{\begin{array}{ll}
    s,\;\, \mbox{if} \;\, G_B=PSL (2,q)_B \rtimes \text{\footnotesize{$\langle$}} \tau_\alpha
\text{\footnotesize{$\rangle$}}\\
    1,\;\,\mbox{if} \;\, G_B=PSL (2,q)_B \rtimes H.\\
\end{array} \right.\]
Clearly, for each $B \in \B$ there exists always a Klein four-group
$V_4 \leq PSL(2,q)$, which fixes $B$ by the definition of Steiner
\mbox{$4$-designs}, and some additional point $x \in X$. We will
distinguish two cases according as $x$ is incident with $B$ or not
and examine for each case the list of possible subgroups of
$PSL(2,q)$ with their orbits on the projective line
(cf.~Lemmas~\ref{PSL_cyc}-\ref{PSL_A_5}). First, let $x \in B$.
Then, clearly \mbox{$k \equiv 1$ (mod $4$)}. It follows that we only
have to consider the subcase when $PSL(2,q)_B$ is conjugate to
$PSL(2,\q)$ with $\q^m = q$, $m \geq 1$. In view of
Lemma~\ref{PSL_kl.PSL}, we obtain $k=\q+1$. Hence,
equation~(\ref{condB}) gives
\begin{equation}\label{eq-0}
(q-2) \left| PSL(2,q)_{0B}\right| = \q(\q-1) (\q -2)s
\end{equation}
\[\mbox{with} \;\, \left| PSL(2,q)_{0B}\right| = \q (\q- 1)\cdot \left\{\begin{array}{ll}
    s,\;\mbox{or}\\
    1.\\
\end{array} \right.\]
Since $q=2^{s^u}$, we can write $\q=2^{s^w}$ for some integer $0
\leq w \leq u$, and $q=\q^m=\q^{s^{u-w}}$. As $k=\q +1 = 2^{s^w}
+1>4$, it follows in particular that $w \geq 1$, and thus $s <
2^{s^w}=\q$. Therefore, using equation~(\ref{eq-0}), we obtain
\[\q^{s^{u-w}}-2=q-2 \leq (\q - 2)s < {\q}^2 -2s.\]
But, as clearly $u-w \geq 1$ (otherwise $k=q+1$, a contradiction to
Corollary~\ref{Cameron_t=4}) this yields a contradiction for every
prime $s$.

Now, let $x \notin B$. Then clearly \mbox{$k \equiv 0$ (mod $4$)}.
We may restrict ourselves to the examination of the following
subcases:

\begin{enumerate}
\item[(i)] $PSL(2,q)_B$ is conjugate to $A_4$ for $s=2$, and $k=12$
in view of Lemma~\ref{PSL_A_4}.

\item[(ii)] $PSL(2,q)_B$ is conjugate to an elementary
Abelian subgroup of order $\q \mid q$ of $PSL(2,q)$, and $k=\q$ due
to Lemma~\ref{PSL_elAb}.

\item[(iii)] $PSL(2,q)_B$ is conjugate to a semi-direct product
of an elementary Abelian subgroup of order $\q \mid q$ with a cyclic
subgroup of order $c$ of $PSL(2,q)$ with $c \mid \q-1$ and $c \mid
q-1$, and $k=\q$ or $\q c$ by Lemma~\ref{PSL_semi}.

\item[(iv)] $PSL(2,q)_B$ is conjugate to $PSL(2,\q)$ with
$\q^m = q$, $m \geq 1$, acting outside the $\q +1$ points mentioned
in the case where $x$ has been incident with $B$, and
Lemma~\ref{PSL_kl.PSL} yields $k=\q(\q-1)$ if $m$ is even, or $k=(\q
+1)\q (\q -1)$.
\end{enumerate}

Again, we can write in the following $\q=2^{s^w}$ for some integer
$0 \leq w \leq u$, and $q=\q^m=\q^{s^{u-w}}$.

(i): Applying equation~(\ref{condB}) yields
\[(q-2) \left| PSL(2,q)_{0B}\right| = 11 \cdot 10 \cdot 9 \cdot 2\]
\[\mbox{with} \;\, \left| PSL(2,q)_{0B}\right| = \left\{\begin{array}{ll}
    2,\;\mbox{or}\\
    1\\
\end{array} \right.,\]
which is clearly impossible.

(iii): Let $k=\q$. By equation~(\ref{condB}), we have
\begin{equation}\label{eq-1}
(q-2) \left| PSL(2,q)_{0B}\right| = (\q-1)(\q-2)(\q-3)  s
\end{equation}
\[\mbox{with} \;\, \left| PSL(2,q)_{0B}\right| = c \cdot \left\{\begin{array}{ll}
    s,\;\mbox{or}\\
    1.\\
\end{array} \right.\]
As $k=\q = 2^{s^w} >4$, we have in particular $w \geq 1$, and hence
$s < 2^{s^w}=\q$. Thus, using equation~(\ref{eq-1}), we obtain
\[q-2 = \q^{s^{u-w}}-2 < \q^3  s < \q^4.\]
Since clearly $u-w \geq 1$ (otherwise $k=q$, which is not possible
by Corollary~\ref{Cameron_t=4}) this yields a contradiction for $s
\geq 5$. If $s=2$, then $\q^{2^{u-w}}-2 < 2\q ^3$ must hold, which
cannot be true for $u-w > 1$. For $s=3$, we may also assume that
$u-w=1$ since otherwise we would have $q=\q^{3^{u-w}} \geq \q^9$,
again a contradiction to the inequality above. As $c \bigm| \q -1$,
it follows for both cases from equation~(\ref{eq-1}) that
\[\q -2 \bigm| q-2,\] and hence
\[2^{s^w-1}-1 \bigm| 2^{s^u-1}-1.\]
Thus, clearly \[s^w-1 \bigm| s^u-1\] and \[w \bigm| u.\] Therefore,
we may conclude that $w=1$ and $u=2$. For $s=2$, it follows that $k=
\q =4$, which has been excluded. For $s=3$, we have $\q = 8$ and
$q=512$, and equation~(\ref{eq-1}) yields
\[510 \cdot \left| PSL(2,q)_{0B}\right| = 7 \cdot 6 \cdot 5 \cdot 3\]
\[\mbox{with} \;\, \left| PSL(2,q)_{0B}\right| =  c  \cdot \left\{\begin{array}{ll}
    3,\;\mbox{or}\\
    1\\
\end{array} \right.,\]
which is clearly impossible.

Now, let $k=\q c$. Then equation~(\ref{condB}) yields
\begin{equation}\label{eq-2}
(q-2) \left| PSL(2,q)_{0B}\right| = (\q c-1)(\q c-2)(\q c-3) s
\end{equation}
\[\mbox{with} \;\, \left| PSL(2,q)_{0B}\right| = \left\{\begin{array}{ll}
    s,\;\mbox{or}\\
    1.\\
\end{array} \right.\]
Polynomial division with remainder gives
\[2^{s^u-1}-1=\bigg(\sum_{i=1}^{\overline{m}}\frac{2^{s^u-1}}{(c \cdot 2^{s^w-1})^i}
\bigg) \bigg(c \cdot 2^{s^w-1}-1\bigg) + \frac{2^{s^u-1}}{(c \cdot
2^{s^w-1})^{\overline{m}}}-1\] for a suitable $\overline{m} \in
\Z^+$ (such that
\[\mbox{deg}\bigg(\frac{2^{s^u-1}}{(c \cdot
2^{s^w-1})^{\overline{m}}}-1 \bigg) < \mbox{deg} \bigg(c \cdot
2^{s^w-1}-1 \bigg)\] as is well-known). As $c$ is odd, clearly
$\big(\frac{2^{s^u-1}}{c \cdot 2^{s^w-1}}\big)^{\overline{m}} \neq
1$, and it follows that $\q c-2$ does not divide $q-2$, yielding a
contradiction to equation~(\ref{eq-2}).

(ii): Let $k=\q$. By equation~(\ref{condB}), we have

\[(q-2) \left|PSL(2,q)_{0B}\right| = (\q-1)(\q-2)(\q-3) s\]
\[\mbox{with} \;\, \left| PSL(2,q)_{0B}\right| = \left\{\begin{array}{ll}
    s,\;\mbox{or}\\
    1.\\
\end{array} \right.\]
As it is easily seen we may argue mutatis mutandis as in
subcase~(iii), $k=\q$.

(iv): If $m>1$ even and $k=\q(\q-1)$, then in view of
equation~(\ref{condB}), we have
\[(q-2) \left| PSL(2,q)_{0B}\right| =
(\q^2-\q-1)(\q^2-\q-2)(\q^2-\q-3)  s\]
\[\mbox{with} \;\, \left| PSL(2,q)_{0B}\right| = (\q+1)\cdot \left\{\begin{array}{ll}
    s,\;\mbox{or}\\
    1.\\
\end{array} \right.\]
As obviously $\q^2-\q-2=(\q+1)(\q-2)$, this leads to a contradiction
analogously as in subcase~(iii), $k=\q$. For $k=\q ^3- \q$,
equation~(\ref{condB}) yields
\begin{equation}\label{eq-3}
(q-2) \left| PSL(2,q)_{0B}\right| = (\q ^3-\q -1)(\q ^3-\q-2)(\q
^3-\q-3) s
\end{equation}
\[\mbox{with} \;\, \left| PSL(2,q)_{0B}\right| = \left\{\begin{array}{ll}
    s,\;\mbox{or}\\
    1.\\
\end{array} \right.\]
We already know that $k=(\q +1)\q(\q-1) \equiv 0$ (mod $4$), and
thus $\q > 2$. If $\left| PSL(2,q)_{0B}\right|=s$, then
\[q = (\q ^3-\q -1)(\q ^3-\q-2)(\q ^3-\q-3)+2= \q ^9- l\]
\[\mbox{with}\;\, l=3\q ^7+6\q ^6-3 \q ^5-12\q^4-10\q^3+6 \q^2 +11 \q+4.\]
As clearly $l>0$, we have $q< \q ^9$. On the other hand, for $\q
> 2$ certainly $l < \q^8(\q-1)$ and hence $q > \q ^8$ must hold, which is impossible
since $q=\q^m$. If $\left| PSL(2,q)_{0B}\right|=1$, then
equation~(\ref{eq-3}) yields
\begin{align*}{}
q= ls+2 \;\,\mbox{with}\;\, l &=(\q ^3-\q -1)(\q ^3-\q-2)(\q
^3-\q-3)\\ &=\q^9-3\q ^7-6\q^6+3 \q ^5+12\q ^4+10\q^3-6 \q^2 -11
\q-6.
\end{align*}
Since $\q = 2^{s^w}>2$, we conclude that $w \geq 1$ and $s <
2^{s^w}=\q$. As obviously $l<\q ^9-1$, it follows  that $q<(\q ^9-1)
\q +2 < \q ^{10}$. On the other hand, for $\q > 2$ clearly $q=ls+2
\geq 2(l+1) > \q ^9$ must hold, which again is impossible.

\pagebreak

\indent \emph{Subcase} $N=PSL(d,q)$, $d \geq 3$.
\smallskip

We have here $\Aut(N)=P \mathit{\Gamma} L(d,q) \rtimes
\text{\footnotesize{$\langle$}} \iota_\beta
\text{\footnotesize{$\rangle$}}$, where $\iota_\beta$ denotes the
graph automorphism induced by the inverse-transpose map
$\beta:GL(d,q) \longrightarrow GL(d,q)$, $x \mapsto {^t(x^{-1})}$.
In the following, let $n=(d,q-1)$.

Let us first assume that $d=3$. In order to show that $G$ with
$PSL(3,q)$ as simple normal subgroup cannot act on any non-trivial
\mbox{$4$-$(q^2+q+1,k,1)$} design, we prove first that $k \leq q+1$.
It is well-known that, for any line $\g$ in the underlying
projective plane $PG(2,q)$, the translation group $T(\g)$ operates
regularly on the points of $PG(2,q)\setminus \g$ and acts trivially
on $\g$. Thus $T(\g)$ fixes a block $B \in \B$ if four or more
distinct points of $B$ lie on $\g$. By the definition of Steiner
\mbox{$4$-designs}, we may choose in $PG(2,q)$ four distinct
collinear points $x_1,x_2,x_3,x_4 \in X$, which are incident with a
unique block $B \in \B$. Let $\g$ denote the line of $PG(2,q)$
through $x_1,x_2,x_3,x_4 \in X$. Consequently, if the block $B$
contains at least one further point of $PG(2,q)\setminus \g$, then
it must contain all points of $PG(2,q)\setminus \g$, thus at least
$q^2+4$ many, which is not possible as $k \leq \bigl\lfloor
\frac{v}{5} + 3 \bigr\rfloor$ by Proposition~\ref{Cam}~(a).
Therefore, $B$ is completely contained in $\g$, and hence  $k \leq
q+1$.

Now, by the definition of Steiner \mbox{$4$-designs}, we may
consider a \mbox{$4$-subset} consisting of three distinct collinear
points $x_1,x_2,x_3 \in X$ and one non-collinear point $x_4 \in X$,
which is incident with a unique block $B \in \B$. If $B$ contains a
fourth point on the line $\g$ of $PG(2,q)$ through $x_1,x_2,x_3 \in
X$, then by the same arguments as above using the translation group
$T(\g)$, we conclude that $B$ lies completely in $\g$, a
contradiction. Thus, we may assume that $B$ contains only further
points which are not on $\g$. Without restriction, we may identify
$x_1=$\text{{$\langle$}}(1,0,0)\text{{$\rangle$}},
$x_2=\text{{$\langle$}}(0,0,1)\text{{$\rangle$}}, x_3 \in
\text{{$\langle$}}x_1,x_2\text{{$\rangle$}},$ and
$x_4=\text{{$\langle$}}(0,1,0)\text{{$\rangle$}}$. As it is known
the cyclic group
\[\left\{ \begin{pmatrix}
  c &        & \\
    & c^{-2} & \\
    &        & c
\end{pmatrix} \Biggm| c \in GF(q)^* \right\} \]
of linear transformations on the associated vector space $V=V(3,q)$
induces a group $U$ of dilatations of order $\frac{q-1}{n}$ on
$PG(2,q)$ with axis the line $\g=\text{{$\langle$}} x_1,x_2
\text{{$\rangle$}}$ and as center the point $x_4$. It is clear that
$U$ fixes each point of its axis as well as its center. Furthermore,
$U$ acts semi-regularly on the points of $PG(2,q) \setminus (\g \cup
\{x_4\})$ and hence all point-orbits on $PG(2,q) \setminus (\g \cup
\{x_4\})$ have length $\frac{q-1}{n}$. As $U$ fixes each of the
points $x_1,x_2,x_3,x_4 \in X$ and hence in particular $B$, we get
\[k \equiv 4 \; \Big(\mbox{mod} \; \frac{q-1}{n}\Big).\]
Due to the fact that $k \leq q+1$, this is obviously impossible if
$3 \nmid q-1$, and for $3 \mid q-1$, we conclude that
\begin{equation} \label{cond}
k=\frac{q-1}{3}+4\;\,\mbox{or}\;\,k=2 \cdot \frac{q-1}{3}+4.
\end{equation}
If we assume that $q>7$, then indeed $q \geq 13$ and we obtain
$\frac{q-1}{3} \geq 4$, which means that we have at least four
distinct collinear points on some line $\h$ of $PG(2,q)$, and we may
argue as above using the translation group $T(\h)$ that $B$ lies
completely in $\h$, which is clearly impossible. Therefore, we only
have to consider the cases when $q=4$ or $7$. For $q=7$,
condition~(\ref{cond}) yields $k=6$ or $8$, whereas $k=6$ can
immediately be ruled out using Lemma~\ref{Comb_t=4}~(d). If any
\mbox{$4$-$(57,8,1)$ design} exists, then there must also exist a
derived \mbox{$3$-$(56,7,1)$ design}. But, for $t=3$, it follows
from Lemma~\ref{Comb_t=4}~(c) that in particular $54$ must be
divisible by $5$, a contradiction. Now, let us assume that $q=4$.
Then only $k=5$ can occur. We have the situation of two intersecting
lines $\g$ and $\h$, and we may distinguish the two cases according
as their intersecting point $x \in \g \cap \h$ is incident with $B$
or not. In the first case, $\g$ and $\h$ are precisely the lines
which intersects $B$ in exactly three distinct points and it can
easily be shown that then $\left| PSL(3,4)_B \right| \leq 8$. In the
second case there is exactly one line which intersects $B$ in
exactly three distinct points and it can be verified here that
$\left| PSL(3,4)_B \right| \leq 2$. However, as there are as blocks
$21$ projective lines in $PG(2,4)$, it follows that $\left|
PSL(3,4)_B \right| \geq \frac{\left|
PSL(3,4)\right|}{b-21}=\frac{20160}{1176}>17$, a contradiction in
both cases. Thus $PSL(3,q)$, and hence also $G$ with $PSL(3,q)$ as
simple normal subgroup, cannot act on any non-trivial
\mbox{$4$-$(q^2+q+1,k,1)$} design.

Now, we consider the case when $d > 3$. Via induction over $d$, we
verify that $G\leq \Aut(\D)$ cannot act on any non-trivial Steiner
\mbox{$4$-design} $\D$. For this, let us assume that there is a
counter-example with $d$ minimal. Without restriction, we can choose
four distinct points $x_1,x_2,x_3,x_4$ from a hyperplane $\h$ of
\mbox{$PG(d-1,q)$}. Analogously as above, it can be shown that the
unique block $B \in \B$ which is incident with the \mbox{$4$-subset}
$\{x_1,x_2,x_3,x_4\}$ is contained completely in $\h$.  Thus, $\h$
induces a $4\mbox{-}(\textstyle{\frac{q^{d-1}-1}{q-1}},k,1) \;
\mbox{design}$, on which $G$ containing $PSL(d-1,q)$ as simple
normal subgroup operates. Inductively, we obtain the minimal
counter-example for $d=3$, a contradiction as above.

\bigskip
\emph{Case} $N=PSU(3,q^2)$, $v=q^3+1$, $q=p^e > 2$.
\medskip

Here $\Aut(N)= P \mathit{\Gamma} U(3,q^2)$, and $\left| G \right| =
(q^3+1)q^3 \frac{(q^2-1)}{n}a$ with $n=(3,q+1)$ and $a \mid 2 ne$.
For the existence of flag-transitive Steiner \mbox{$4$-designs},
necessarily
\[ r = \frac{q^3(q^3-1)(q^3-2)}{(k-1)(k-2)(k-3)} \Bigm|
\left| G_0 \right| \Bigm| \left| P \mathit{\Gamma} U(3,q^2)_0
\right| = q^3(q^2-1)2e \] must hold in view of Lemma~\ref{divprop}.
As obviously $(q^2+q+1,q+1)=1$ and $(q^3-2,q+1)=(3,q+1)=n$, we have
in particular
\begin{equation}\label{eq_PSU}
(q^3-2) (q^2+q+1)\bigm| (k-1)(k-2)(k-3)2ne,\; \mbox{where}\,\; e
\leq \mbox{log}_2q.
\end{equation}
On the other hand, Corollary~\ref{Cameron_t=4} yields $k \leq
\bigl\lfloor \sqrt{q^3+1} + \frac{5}{2}\bigr\rfloor <
q^{\frac{3}{2}}+3$. Hence, using property~(\ref{eq_PSU}), we have
only a small number of possibilities to check, which can easily be
eliminated by hand.

\bigskip
\emph{Case} $N=Sz(q)$, $v=q^2+1$, $q=2^{2e+1}>2$.
\medskip

We have $\Aut(N)= Sz(q) \rtimes \text{\footnotesize{$\langle$}}
\alpha \text{\footnotesize{$\rangle$}}$, where $\alpha$ denotes the
Frobenius automorphism $GF(q) \longrightarrow GF(q)$, \mbox{$x
\mapsto x^2$}. Thus, by Dedekind's law, $G = Sz(q) \rtimes (G \cap
\text{\footnotesize{$\langle$}} \alpha
\text{\footnotesize{$\rangle$}})$, and $\left| G \right| =
(q^2+1)q^2(q-1)a$ with $a \mid 2e+1$. From Remark~\ref{equa_t=4}, we
obtain
\[(q^2-2)(q+1) = (k-1)(k-2)(k-3) \frac{a}{\left|G_{xB}\right|}
\;\, \mbox{if} \;\, x \in B.\]

First, we show that every element $g \in G$ that fixes three
distinct points must fix at least five distinct points. Let us
assume that $g \in G$ with $\left| \mbox{Fix}_X(g)\right| \geq 3$.
Let $x \in \mbox{Fix}_X(g)$, and $P$ the normal Sylow
\mbox{$2$-subgroup} of $Sz(q)_x$ acting regularly on $X\setminus
\{x\}$. Set $C:=C_P (g)$. If \mbox{$y,z \in \mbox{Fix}_X(g)
\setminus \{x\}$}, then $z=y^h$ with $h \in P$. Thus, as
$y^{hg}=y^h=y^{gh}$, we conclude that
\[[h^{-1}, g^{-1}] \in G_{xy} \cap [P, G_x] \leq P_y=1.\]
Then $h\in C$, and hence $C$ acts point-transitively on
\mbox{$\mbox{Fix}_X(g) \setminus \{x\}$}. Therefore, as $\left|
\mbox{Fix}_X(g)\right| \geq 3$, it follows that
\mbox{$\left|\mbox{Fix}_X(g)\right| \equiv 1$ (mod $2$).} Clearly,
the set $\mbox{Fix}_X(g)$ is left invariant by $C_{Sz(q)}(g)$ and
$C_{Sz(q)}(g)$ operates on $\mbox{Fix}_X(g)$. Since $x \in
\mbox{Fix}_X(g)$ can be chosen arbitrarily, it follows that
$C_{Sz(q)}(g)$ operates point-transitively on $\mbox{Fix}_X(g)$, and
thus $\left| \mbox{Fix}_X(g)\right| \big| \left| Sz(q)\right|$. As
the order of $Sz(q)$ is not divisible by $3$, clearly $\left|
\mbox{Fix}_X(g)\right| \neq 3$, and due to the fact that
\mbox{$\left|\mbox{Fix}_X(g)\right| \equiv 1$ (mod $2$),} we have
$\left|\mbox{Fix}_X(g)\right| \geq 5$.

Since $G$ is block-transitive, it is sufficient to consider some
appropriate, unique block $B \in \B$. As clearly
$\text{\footnotesize{$\langle$}} \alpha
\text{\footnotesize{$\rangle$}} \leq \Aut(N)_{0,1,\infty}$, it
follows from above that $\text{\footnotesize{$\langle$}} \alpha
\text{\footnotesize{$\rangle$}}$ must fix some fourth point, and
hence $G \cap \text{\footnotesize{$\langle$}} \alpha
\text{\footnotesize{$\rangle$}} \leq G_{0B}$ by the definition of
Steiner \mbox{$4$-designs}. Thus, we have particularly
\[(q^2-2)(q+1) \leq (k-1)(k-2)(k-3),\]
which does obviously not hold for $k \leq q+2$. On the other hand,
Corollary~\ref{Cameron_t=4} yields $k \leq \bigl\lfloor
\sqrt{q^2+1}+\frac{5}{2} \bigr\rfloor < q+3$, a contradiction.

\bigskip
\emph{Case} $N=Re(q)$, $v=q^3+1$, $q=3^{2e+1}>3$.
\medskip

Here $\Aut(N)= Re(q) \rtimes \text{\footnotesize{$\langle$}} \alpha
\text{\footnotesize{$\rangle$}}$, where $\alpha$ denotes the
Frobenius automorphism $GF(q) \longrightarrow GF(q)$, \mbox{$x
\mapsto x^3$}. Thus, by Dedekind's law, $G = Re(q) \rtimes (G \cap
\text{\footnotesize{$\langle$}} \alpha
\text{\footnotesize{$\rangle$}})$, and $\left| G \right| =
(q^3+1)q^3(q-1)a$ with $a \mid 2e+1$. It follows from
Remark~\ref{equa_t=4} that
\[(q^3-2)(q^2+q+1) = (k-1)(k-2)(k-3) \frac{a}{\left|G_{xB}\right|}
\;\, \mbox{if} \;\, x \in B.\] First, we show that every element $g
\in G$ that fixes three distinct points must also fix a fourth
point. Let us assume that $g \in G$ with $\left|
\mbox{Fix}_X(g)\right| \geq 3$. Let $x \in \mbox{Fix}_X(g)$, and $P$
the normal Sylow \mbox{$3$-subgroup} of $Re(q)_x$ acting regularly
on $X\setminus \{x\}$. As in the previous Case, it can be shown that
then $C_P (g)$ acts point-transitively on \mbox{$\mbox{Fix}_X(g)
\setminus \{x\}$}. Thus, we have \mbox{$\left|\mbox{Fix}_X(g)\right|
\equiv 0$ (mod $2$)}, and the claim follows.

Since $G$ is block-transitive, it is sufficient to consider some
appropriate, unique block $B \in \B$. As clearly
$\text{\footnotesize{$\langle$}} \alpha
\text{\footnotesize{$\rangle$}} \leq \Aut(N)_{0,1,\infty}$, we
deduce from above that $G \cap \text{\footnotesize{$\langle$}}
\alpha \text{\footnotesize{$\rangle$}} \leq G_{0B}$ by the
definition of Steiner \mbox{$4$-designs}. Hence, we have in
particular
\[(q^3-2)(q^2+q+1) \leq (k-1)(k-2)(k-3),\]
which is not possible as Corollary~\ref{Cameron_t=4} yields $k \leq
\bigl\lfloor \sqrt{q^3+1} + \frac{5}{2} \bigr\rfloor <
q^{\frac{3}{2}}+3$.

\bigskip
\emph{Case} $N=Sp(2d,2)$, $d \geq 3$, $v = 2^{2d-1} \pm 2^{d-1}$.
\medskip

As here \mbox{$\left|\Out(N) \right|=1$}, we have $N=G$. Let $X^+$
respectively $X^-$ denote the set of points on which $G$ operates.
It is well-known that $G_z$ acts on \linebreak $X^{\pm} \setminus
\{z\}$ as $O^\pm(2d,2)$ does in its usual rank $3$ manner on
singular points of the underlying non-degenerate orthogonal space
$V^\pm = V^\pm (2d,2)$.

It is easily seen that there are $2^{2d-2}(2^d \mp 1)(2^{d-1} \pm
1)$ hyperbolic pairs in $V^\pm$, and by Witt's theorem, $O^\pm
(2d,2)$ is transitive on these hyperbolic pairs. Let $\{x,y\}$
denote a hyperbolic pair, and
$\mathcal{E}=\text{\footnotesize{$\langle$}} x,y
\text{\footnotesize{$\rangle$}}$ the hyperbolic plane spanned by
$\{x,y\}$. As $\mathcal{E}$ is non-degenerate, we have the
orthogonal decomposition
\[V^\pm=\mathcal{E}\perp\mathcal{E}^\perp.\]
Clearly, $O^\pm(2d,2)_{\{x,y\}}$ stabilizes $\mathcal{E}^\perp$
as a subspace, which implies that
\linebreak \mbox{$O^\pm(2d,2)_{\{x,y\}} \cong O^\pm(2d-2,2)$}.
Therefore, we have
\[O^\pm(2d-2,2)\cong O^\pm(2d,2)_{\{x,y\}} \unlhd O^\pm(2d,2)_\mathcal{E} = G_{z,\mathcal{E}}.\]
Since $O^\pm(2d-2,2)$ acts transitively on the singular points of
the $(2d-2)$-dimensional orthogonal subspace, the smallest orbit on
$V^\pm \setminus \mathcal{E}$ under $G_{z,\mathcal{E}}$ has length
at least $2^{2d-3} \pm 2^{d-2}$. If the unique block $B \in \B$
which is incident with the \mbox{$4$-subset} \mbox{$\{x,y,x+y,z\}$}
contains some singular point in \mbox{$V^\pm \setminus
\mathcal{E}$}, then we would have $k \geq 2^{2d-3} \pm 2^{d-2}+4$, a
contradiction to Corollary~\ref{Cameron_t=4}. Thus, all points of
$B$ apart from $z$ lie completely in $\mathcal{E}$. By the
flag-transitivity of $G$, it follows that for each block all points
apart from a singleton must be contained in an affine plane. Thus
$k=5$, which is impossible since \mbox{$k \equiv 0$ (mod $4$)} by
Lemma~\ref{Comb_t=4}~(c).

\bigskip
\emph{Cases} $N=PSL(2,11)$, $v=11$, and  $N=M_{11}$, $v=12$.
\medskip

As is known, these exceptional permutation actions occur inside the
Mathieu group $M_{24}$ in its action on $24$ points. This set can be
partitioned into two sets $X_1$ and $X_2$ of $12$ points each such
that the setwise stabilizer of $X_1$ is the Mathieu group $M_{12}$.
The stabilizer in this latter group of a point $x$ in $X_1$ is
isomorphic to $M_{11}$ and operates (apart from its natural
\mbox{$4$-transitive} action on \mbox{$X_1 \setminus \{x\}$})
\mbox{$3$-transitively} on the $12$ points of $X_2$. The one-point
stabilizer in this action of degree $12$ is $PSL(2,11)$ acting
$2$-transitively on $11$ points. The geometry preserved by the
$3$-transitive action of $M_{11}$ is not a Steiner
\mbox{$t$-design}, but a $3$-$(12,6,2)$ design
(e.g.~\cite[Ch.\,IV,\,5.3]{BJL1999}). Hence also the derived design
$\D$ on which \mbox{$G \leq \Aut(\D)$} acts cannot be a Steiner
design.

\pagebreak

\emph{Case} $N=M_v$, $v=11,12,22,23,24$.
\medskip

If $v=11,12,23$ or $24$, then $G=M_v$ is always
\mbox{$4$-transitive} and thus~\cite[Thm.\,3]{Kant1985} yields the
designs described in the Main Theorem. Obviously, flag-transitivity
holds as the $4$-transitivity of $G$ implies that $G_x$ acts
block-transitively on the derived Steiner \mbox{$3$-design} $\D_x$
for any $x \in X$. For $v=22$, Corollary~\ref{Cameron_t=4} gives $k
\leq 7$, and again the cases for $k$ can easily be eliminated in
view of Lemma~\ref{Comb_t=4}~(c) and~(d).

\medskip

This completes the proof of the Main Theorem.

\medskip


\subsection*{Acknowledgment}
\mbox{I am grateful to C. Hering for helpful conversations.}
\bibliographystyle{amsplain}
\bibliography{Xbib4des}
\end{document}